\documentclass[11pt,reqno]{amsart}
\usepackage{amsmath,amssymb,latexsym}
\usepackage[small]{caption}
\usepackage{graphicx,color,mathrsfs,tikz}\usepackage{subfigure,color}

\usepackage{cite}
\usepackage[colorlinks=true,urlcolor=blue,
citecolor=red,linkcolor=blue,linktocpage,pdfpagelabels,
bookmarksnumbered,bookmarksopen]{hyperref}
\usepackage{units}
\usepackage{enumitem}
\usepackage[left=2.6cm,right=2.6cm,top=2.75cm,bottom=2.75cm]{geometry}
\usepackage[hyperpageref]{backref}

\usepackage[colorinlistoftodos,prependcaption]{todonotes}

\numberwithin{equation}{section}
\newtheorem{theorem}{Theorem}[section]
\newtheorem{proposition}[theorem]{Proposition}
\newtheorem{lemma}[theorem]{Lemma}
\newtheorem{remark}[theorem]{Remark}

\newtheorem{corollary}[theorem]{Corollary}
\newtheorem{definition}[theorem]{Definition}
\theoremstyle{definition}

\newtheorem{claim}{Claim}

\newcommand{\I}{\infty}
\renewcommand{\S}{\sigma}

\newcommand{\T}{\tilde}

\newcommand{\N}{\nabla}

\renewcommand{\L}{\lambda}

\renewcommand{\r}{\Lambda}

\renewcommand{\P}{\partial}
\newcommand{\B}{\beta}

\newcommand{\Lap}{\triangle}

\newcommand{\A}{\alpha}
\newcommand{\lng}{\langle}
\newcommand{\rng}{\rangle}

\newcommand{\V}{\varepsilon}
\newcommand{\vp}{\varphi}

\newcommand{\QED}{\qquad\mbox{\qed}}

\newcommand{\R}{\mathbb{R}}
\newcommand{\BN}{\mathbb{N}}

\newcommand{\CB}{\mathcal{B}}
\newcommand{\CG}{\mathcal{G}}

\newcommand{\CF}{\mathcal{F}}
\newcommand{\CL}{\mathcal{L}}

\newcommand{\CJ}{\mathcal{J}}
\newcommand{\CK}{\mathcal{K}}
\newcommand{\CP}{\mathcal{P}}
\newcommand{\CR}{\mathcal{R}}

\newcommand{\CM}{\mathcal{M}}

\renewcommand{\CG}{\mathcal{G}}

\newcommand{\CZ}{\mathcal{Z}}

\DeclareMathOperator{\supp}{supp}

\title[Weighted estimates and Fujita exponent]{Weighted estimates and Fujita exponent for a nonlocal equation}

\author[S.~Khomrutai]{Sujin Khomrutai}

\address[S.~Khomrutai]{Department of Mathematics and Computer Science,\newline\indent
 Faculty of Science, Chulalongkorn  University, \newline\indent
 Bangkok 10330, Thailand
}
\email{sujin.k@chula.ac.th}



\begin{document}

\begin{abstract}
We investigate a nonlocal equation $\partial_tu=\int_{\R^n}J(x-y)(u(y,t)-u(x,t))dy+a(x,t)u^p$ in $\mathbb{R}^n$, where $a$ is unbounded and $J$ belongs to a weighted space. Crucial
weighted $L^p$ and interpolation estimates for the Green operator are established by a new method based on the sharp Young's inequality, the asymptotic behavior of a regular varying coefficients exponential series, and the properties of auxiliary functions $\varGamma=(1+|x|^2/\eta)^{b/2}$ that $-\varGamma/\eta\lesssim J\ast\varGamma-\varGamma\lesssim\varGamma/\eta$ and $\eta^{-b_+/2}\lesssim\varGamma/\left\langle x\right\rangle^b\lesssim \eta^{-b_-/2}$. 
Blow-up behaviors are investigated by employing test functions $\phi_R=\varGamma$ ($\eta=R$) instead of principal eigenfunctions. Global well-posedness in weighted $L^p$ spaces for the Cauchy problem is proved.  When $a\sim\left\langle x\right\rangle^\sigma$ the Fujita exponent is shown to be $1+(\sigma+2)/n$. Our approach generalizes and unifies nonlocal diffusion equations and pseudoparabolic equations.
\end{abstract}
\maketitle

\section{Introduction}
\subsection{Overview}
Semilinear equations of the form
\begin{align}\label{Eqn:GenL}
\P_tu=\CL u+a(x,t)u^p\quad(x\in\R^n,t>0),
\end{align}
where $\CL$ is a differential or nonlocal operator,
have been the subject of many recent articles on nonlinear evolution equations. For the semilinear heat equation, the result is now classic and it is known that the Fujita exponent is $1+2/n$ when $a=1$, and $1+(\sigma+2)/n$ when $a\sim\lng x\rng^\sigma$. See \cite{Fujita66, Hayakawa73, KobayashiEtal77, AronsonEtal78, Weissler81} and \cite{Levine90,Pinsky97,DengEtal00,Suzuki02,
FerreiraEtal06, QuittnerEtal07, BaiEtal11} for more details and some generalizations. 
For the fractional heat equation $\CL=-(-\Lap)^s$ ($0<s<1$) with $a=1$, it was found \cite{Sugitani75, Birkner02} that the Fujita exponent is $1+2s/n$. Self-similarity and Fourier transform were used to obtain these results.\medskip

In this work, we study a nonlocal equation of the form (\ref{Eqn:GenL}), i.e.
\[
\CL u=\int_{\R^n}J(x-y)(u(y)-u(x))dx,
\]
where $J$ is an integrable function. We generalize many recent results in the literature. 
In \cite{MelianEtal10}, the Fujita exponent of (\ref{Eqn:GenL}) was shown to be $1+2/n$ when $a=1$ and
$J=J(|x|)\geq0$ is a continuous, decreasing, and compactly supported function. 
Kaplan's eigenfunction method \cite{Kaplan63} (with results from \cite{MelianEtal09}) and 
the test function method \cite{MitidieriEtal01} were employed in the blow-up analysis of \cite{MelianEtal10}. We note that this blow-up analysis (and results in \cite{MelianEtal09}) relies strongly on the condition of $J$.
Global well-posedness in $L^p$ spaces was obtained by a
comparison argument.
In \cite{YangEtal14}, assume that $a=a(x)$ and $J=J(|x|)\geq0$ are compactly supported, continuous functions, it was shown that (\ref{Eqn:GenL}) exhibits 
blow-up in finite time precisely when $n=1$. Finally, we mention the recent work \cite{Alfaro17} when $a=1$ and the Fourier transform of $J$ satisfies
\begin{align}\label{Four:Jhat}
\widehat{J}(\xi)=1-A|\xi|^{2s}+o(|\xi|^{2s}),\quad0<s\leq1,\,\,\mbox{as $|\xi|\to0$}.
\end{align}
The Fujita exponent in this case was found to be $1+2s/n$, same as the fractional heat equation.\medskip

In addition to the preceding nonlocal equations, our results also generalize some weakly singular equations and higher order PDEs that can be put into a nonlocal form. The third order semilinear pseudoparabolic equation studied in \cite{CaoYinWang09,KhomrutaiNA15}
\begin{align}\label{Eqn:pseudo}
\P_tu-\Lap\P_tu=\Lap u+a(x,t)u^p\quad(x\in\R^n,t>0)
\end{align}
can be expressed in a nonlocal form as
\begin{align*}
\P_tu=\int_{\R^n}B(x-y)(u(y,t)-u(x,t))dy+B\ast(au^p),
\end{align*}
where $B$ is the kernel of the Bessel potential operator $(1-\Lap)^{-1}$. Assuming $a\sim\lng x\rng^\sigma$ ($\sigma\geq0$), it was shown \cite{KhomrutaiNA15} (see also\cite{CaoYinWang09} when $a=1$) that the Fujita exponent for (\ref{Eqn:pseudo}) is $1+(\sigma+2)/n$. Since 
\begin{align}\label{Four:B}
B\ast\vp
=\CF^{-1}((1+|\xi|^2)^{-1}\widehat{\vp}),
\end{align}
weighted $L^p$ and interpolation estimates, which are crucial 
for the local and global well-posedness results, were proved in \cite{KhomrutaiNA15} (see also \cite{HayashiEtal06} when $a=1$) by Fourier transform techniques. 
We shall, in this work, reprove these estimates from the other perspective which depends only upon the singularity and behavior as $|x|\to\infty$ of the kernel $B$. This will open the way up for the study of pseudoparabolic equation with a more general (variable coefficients) elliptic operator rather than the Laplacian.\medskip

Now let us state the main problem of this work. We study the semilinear nonlocal equation
\begin{align}\tag{NL}\label{Eqn:main}
\begin{cases}
\displaystyle
\P_tu=\int_{\R^n}J(x-y)(u(y,t)-u(x,t))dy+a(x,t)u^p\quad (x\in\R^n,t>0),\\
\vspace{-10pt}\\
\displaystyle
u(x,0)=u_0(x),\quad x\in\R^n.
\end{cases}
\end{align}
The inhomogeneous term can be replaced by $\CK\ast(au^p)$, where $\CK$ is a bounded linear operator. 
Evolution equations of this type have been found to be connected with many physical phenomena, see for instance \cite{AndreuEtal10,AKS11} and the references therein. 
Unlike \cite{MelianEtal10,TerraEtal11,
YangEtal14,Yang15}, the kernel $J$ in this work needs not be compactly supported nor radially decreasing, and it can be weakly singular. 
Moreover, the inhomogeneous term coefficient $a$ can be unbounded. Under this circumstance, we are facing some difficulties. 
As in \cite{KhomrutaiNA15,HayashiEtal06}, we work around the unboundedness of $a$ by investigating (\ref{Eqn:main}) within weighted $L^p$ spaces. Moreover, for such a general $J$, the limit properties of principal eigenvalues and eigenfunctions on dilated domains, that are needed in conjunction with Kaplan's method for the blow-up analysis, are unavailable. In fact, the properties obtained in\cite{MelianEtal09,Coville10} require $J$ to be compactly supported and radially decreasing. We get around this difficulty by employing test functions $\phi_R=(1+|x|^2/R)^{-b/2}$ on $\R^n$ rather than principal eigenfunctions. Finally, in this work, there is no property like (\ref{Four:Jhat}) or (\ref{Four:B}) for the kernel, which prevents us from using Fourier transform tools, and no obvious self-similarity property for (\ref{Eqn:main}). We shall prove the crucial weighted $L^p$ and interpolation estimates via a purely analytical argument.

\subsection{Main results}

We investigate mild solutions for (\ref{Eqn:main}). For the homogeneous equation, the ``diffusive" kernel $J$ leads to the solution of the initial value problem
\begin{align*}
\P_tu=\CL u,\quad u|_{t=0}=u_0,
\end{align*}
expressed in terms of the Green operator as
\begin{align}\label{GreenDef}
\CG(t)u_0=e^{(\CJ-\A_0)t}:=e^{-\A_0t}\sum_{k=0}^\infty\frac{t^k}{k!}\CJ^ku_0,
\end{align}
where 
\[
\A_0=\|J\|_{L^1},\quad\CJ\vp=\int_{\R^n}J(x-y)\vp(y)dy,\quad\CJ^0=id,\quad \CJ^k=J^{\ast k}\,\,(k\geq1).
\]

We first prove the following weighted $L^p$ estimate.

\begin{theorem}\label{Thm:GreenFar}

Assume that there exists $\delta\geq2$ such that
\begin{align}\label{Est:Jbdd}
J\in L^{1}_{\delta}(\R^n),\,\,\mbox{i.e.}\quad
\int_{\R^n}|J(x)|\lng x\rng^\delta dx<\infty.
\end{align}
Then, for any $b\in\R$ with $|b|\leq\delta-2$ and $1\leq q\leq\infty$, we have 
\begin{align}\label{Est:ThmGreenFar1}
\|\CG(t)f\|_{L^q_b}\leq C
\lng t\rng^{|b|/2}\|f\|_{L^q_b}.
\end{align}

\end{theorem}

Then we immediately obtain the following local well-posedness result.

\begin{theorem}\label{Thm:LocWell}

Let $p>1$. Assume $J\in L^1_\delta(\R^n)$ with $\delta\geq2$ and $a(x,t)$ satisfies
\begin{align}\label{Hyp:LocH1}
|a(x,t)|\lesssim\lng x\rng^\sigma\quad(\forall\,x\in\R^n,t>0),
\end{align}
where $\sigma\in\R$ and
\begin{align}\label{Hyp:LocH2}
\frac{\sigma_+}{p-1}\leq\delta-2.
\end{align}
Assume $u_0\in L^\infty_b(\R^n)\cap C(\R^n)$ where $b\in\R$ is such that 
\begin{align}\label{Hyp:LocH3}
|b|\leq\delta-2\quad\mbox{and}\quad b\geq\frac{\sigma}{p-1}.
\end{align}
Then there exists $T>0$ such that (\ref{Eqn:main}) has a unique solution $u\in C\left([0,T];L^\infty_b(\R^n)\cap C(\R^n)\right)$. If in addition $J$ and $a(x,t)$ are non-negative functions, then $u(x,t)\geq0$.

\end{theorem}

Theorem \ref{Thm:LocWell} is proved by a standard Banach contraction mapping argument using the bound (\ref{Est:ThmGreenFar1}) of $\CG(t)$. 
To prove Theorem \ref{Thm:GreenFar}, we show in Theorem \ref{Thm:AsympPsi} that auxiliary functions
\[
\varGamma=\left(1+\frac{|x|^2}{\eta}\right)^{b/2}\quad(\eta\geq2,b\in\R),
\]
satisfy
\begin{align}\label{Ineq:LT}
-\frac{d}{\eta}\varGamma\leq J\ast\varGamma-\A_0\varGamma\leq\frac{d}{\eta}\varGamma,
\end{align}
where $d=d(b,J,\A_0)\geq0$ is a constant. These inequalities can be interpreted as that the convolution with $J$ increases or decreases $\varGamma$ by a factor in $[\A_0-d/\eta,\A_0+d/\eta]$. So, in particular, as $\eta\to\infty$, $\varGamma$ approaches a stationary solution of $\P_tu=\CL u$. 
The choice of auxiliary functions to be used in proving Theorem \ref{Thm:GreenFar} (and Theorem \ref{Thm:Interp})
comes from the following observation. If $a\sim\lng x\rng^b={}$the weight, then $\varGamma\sim a$ according to the inequalities
\begin{align}\label{Ineq:Tweight}
\eta^{-b_+/2}\leq\frac{\varGamma}{\lng x\rng^b}\leq \eta^{-b_-/2}.
\end{align}
See Lemma \ref{Lem:GBr}.\medskip

Applying (\ref{Ineq:LT}) and (\ref{Ineq:Tweight}) with $\eta=\eta_0+t$, we can prove (\ref{Est:ThmGreenFar1}) for $q=1$ and $q=\infty$ separately. The Riesz-Thorin interpolation theorem then implies the desired estimate for any $1\leq q\leq\infty$. 

\begin{remark}\rm 
Unless $J$ satisfies a certain stronger property (see Theorem \ref{Thm:Interp}), (\ref{Est:ThmGreenFar1}) is optimal, at least when $q=\infty$, as exhibited in unbounded solutions in Theorem 1.4 of \cite{BrandleEtal11}.
If $b=0$, it also implies a decay in $L^q$-norm of the Green operator. 
See Corollary 2.6 in \cite{CordobaEtal04}, for the decay in $L^q$-norm estimate for the fractional heat equation.
  
\end{remark}

Our next main result is the following weighted interpolation estimate. 

\begin{theorem}\label{Thm:Interp}

Assume there exist $\varepsilon_0>0$ and $\beta\in(n,\infty]$ such that
\begin{align}\label{Est:InterpHypJ}
J\in L^1_{2+\beta}(\R^n)\cap L^{1+\varepsilon_0}_\beta(\R^n),
\end{align}
i.e.
\[
\max\left\{\int_{\R^n}|J(x)|\lng x\rng^{2+\beta}dx,\int_{\R^n}(|J(x)|\lng x\rng^{\beta})^{1+\varepsilon_0}dx\right\}<\infty.
\]
Let $1\leq q\leq Q\leq\infty$ and $b\in\R$ with
\begin{align}
|b|<\beta-n\left(\frac{1}{Q}-\frac{1}{q}+1\right).
\end{align}
Then we have
\begin{align}\label{Est:IntqQ}
\|\CG(t)f\|_{L^{Q}_b}\lesssim\lng t\rng^{\frac{n}{2}\left(\frac{1}{Q}-\frac{1}{q}+\frac{|b|}{n}\right)}\|f\|_{L^q}+\lng t\rng^{\frac{n}{2}\left(\frac{1}{Q}-\frac{1}{q}\right)}\|f\|_{L^q_b}+e^{-t/2}\|f\|_{L^{Q}_b}.
\end{align}

\end{theorem}

Then we can prove the following global well-posedness result.

\begin{theorem}\label{Thm:Global}

Assume $J\geq0$ and there exists $\varepsilon_0>0$ such that
\begin{align}\label{Hyp:JGlobal}
J\in L^1_\infty(\R^n)\cap L^{1+\varepsilon_0}_\infty(\R^n),
\end{align}
i.e.~$J\in L^1_\beta(\R^n)\cap L^{1+\varepsilon_0}_\beta(\R^n)$ for all $\beta>0$, and there is $\sigma\geq0$ such that
\begin{align}\label{Hyp:global}
a(x,t)\sim\lng x\rng^\sigma.
\end{align}
Let $p>1+\frac{\S+2}{n}$. 
If $u_0\geq0$ and $\|u_0\|_{L^1_{\S/(p-1)}}+\|u_0\|_{L^\infty_{\S/(p-1)}}$ is sufficiently small, then (\ref{Eqn:main}) has a unique global solution 
\[
u\in C\left([0,\infty);C(\R^n)\cap L^1_{\S/(p-1)}(\R^n)\cap L^\infty_{\S/(p-1)}(\R^n)\right).
\]

\end{theorem}
The procedure to prove Theorem \ref{Thm:Global} is more or less standard. We choose a suitable mixed-norm involving the $L^1_{\sigma/(p-1)}$- and $L^\infty_{\sigma/(p-1)}$-norms and apply Theorem \ref{Thm:Interp} to bound the various norms. See for instance \cite{HayashiEtal06,KhomrutaiNA15,QuittnerEtal07}.\medskip

To prove Theorem \ref{Thm:Interp}, we choose $N\geq1/\varepsilon_0+1$ and split $\CG(t)=\CG_N(t)+\CR_N(t)$, where $\CG_N$ has the first $N$ terms and $\CR_N$ is the remainder. The key for proving the interpolation estimate (\ref{Est:IntqQ}) is a precise pointwise decay estimate for the kernel $R_N(x,t)$ of $\CR_N$, which exists for such a large $N$. In fact, we prove in Proposition \ref{Prop:RN} that
\begin{align}\label{RNdecay}
|R_N(x,t)|\lesssim \vartheta^{-\beta/2}t^{-n/2}\quad\mbox{as $|x|,t,\vartheta\geq1$},
\end{align}
where $\vartheta:=|x|^2/t$ is the self-similar variable.
Then, using this decay estimate, we can prove (\ref{Est:IntqQ}) directly. To the best of the author knowledge, the above pointwise decay estimate of $R_N$ and the weighted interpolation estimate (\ref{Est:IntqQ}) for nonlocal operators are new.

\begin{remark}\rm
A decay estimate similar to (\ref{RNdecay}) associated with a differential or pseudodifferential operator is often derived by Fourier transform technique (see \cite{HayashiEtal06,KhomrutaiNA15}). In this work, we derive (\ref{RNdecay}) by applying the sharp Young's inequality (Proposition \ref{Prop:SharpYoung}), the inequality (\ref{Ineq:LT}), and the bound of regular varying coefficients exponential series (Proposition \ref{Prop:Kummer}). 
The hypothesis (\ref{Hyp:JGlobal}) can be weaken a bit but we choose not do to simplify the presentation.

\end{remark}

\begin{remark}\rm

The hypotheses of Theorem \ref{Thm:GreenFar} and Theorem \ref{Thm:Interp} are satisfied by the Bessel potential kernel $B$ discussed in (\ref{Eqn:pseudo}). In fact, $B\in L^1_\infty(\R^n)\cap L^{1+\varepsilon_0}_\infty(\R^n)$ for some $\varepsilon_0>0$. So technique in this work can be extended to the pseudoparabilic equation (\ref{Eqn:pseudo}) and many other problems, for instance, pseudoparabolic regularization of parabolic equation $\P_tu-\V\CP\P_tu=\CP u+f$, where $\CP u=\sum_{i,j}a_{ij}(x)\P_{ij}u+\sum_ib_i(x)\P_iu+c(x)u$ is a second order elliptic operator with non-constant coefficients.

\end{remark}

Next we investigate the blow up behaviors for solutions of (\ref{Eqn:main}).

\begin{theorem}\label{Thm:Blow-up}
Assume $J\geq0$, $J\in L^1_\infty(\R^n)$, and the coefficient $a(x,t)$ satisfies
\begin{align}\label{Hyp:blowup}
a(x,t)\geq C\lng x\rng^\sigma,\quad\sigma>-2.
\end{align}
\begin{itemize}
\item[(1)] (sub-critical exponent) Let $1<p<1+(\sigma+2)/n$. If $u_0\geq0$ is nontrivial, then the solution $u$ of (\ref{Eqn:main}) blows up in finite time in both $L^1$- and $L^\infty$-norms, i.e.\  $\exists\,T_0<\infty$ such that
\begin{align}
\|u(t)\|_{L^1}\to\infty\quad\mbox{and}\quad\|u(t)\|_{L^\infty}\to\infty,\quad\mbox{as $t\to T_0^-$}.
\end{align}
\item[(2)] (critical exponent) Let $p=1+(\sigma+2)/n$. If $u_0\geq0$ is nontrivial, then the solution $u$ of (\ref{Eqn:main}) blows up in finite time in both $L^1$- and $L^\infty$-norms.
\item[(3)] (super-critical exponent) Let $p>1+(\sigma+2)/n$. There are constants $b_0=b_0(n,\S,p)>0$ and $m_0=m_0(n,b,\S,p)>0$ such that if $u_0\geq0$ satisfies
\begin{align}
\int_{\R^n}\left(1+\frac{|x|^2}{2}\right)^{-b}u_0(x)dx>m_0
\end{align}
for some $b>b_0$, then the solution $u$ blows up in finite time.
\end{itemize}

\end{theorem}

We prove Theorem \ref{Thm:Blow-up} as follows. Define 
\begin{align}
\phi_R(x)=\left(1+\frac{|x|^2}{R}\right)^{-b/2}\,\,(b>n)\quad\mbox{and}\quad f_R(t):=\int_{\R^n}\phi_Rudx.
\end{align}
Multiplying (\ref{Eqn:main}) with the test function $\phi_R$, integrating both sides of the equation over $\R^n$, 
and applying Lemma \ref{Lem:Blowup1}, 
H\"older's inequality, we obtain the Bernoulli type differential inequality
\begin{align}\label{Est:Bern}
\frac{d}{dt}f_R\geq-\L_Rf_R+C\mu_Rf_R^p,
\end{align}
where
\[
\L_R:=\frac{d}{R}\quad\mbox{and}\quad\mu_R:=\begin{cases}
1&\mbox{if $1<p<1+\frac{\sigma}{n}$},\\
(\ln R)^{1-p}&\mbox{if $p=1+\frac{\sigma}{n}$},\\
R^{-\frac{1}{2}(n(p-1)-\sigma)}&\mbox{if $p>1+\frac{\sigma}{n}$}.
\end{cases}
\]

If $1<p<p_F:=1+(\sigma+2)/n$, then we have $\L_R/\mu_R\to0$ as $R\to\infty$. So by taking $R$ sufficiently large, it follows from Lemma \ref{Lem:Blowup2} and (\ref{Est:Bern}) that $f_R(t)\to\infty$ at a finite time $T_0$. Using that
$f_R(t)\leq\max\{\|u(t)\|_{L^1},\|\phi_R\|_{L^1}\|u(t)\|_{L^\infty}\}$, 
the blow-up in both $L^1$- and $L^\infty$-norms is true, hence proving Theorem \ref{Thm:Blow-up} (1).\medskip

The argument to prove Theorem \ref{Thm:Blow-up} (2) is a bit more delicate. It starts by noticing that for $p=p_F$, the differential inequality (\ref{Est:Bern}) has $\mu_R=1/R$. 
So the blow-up criterion obtained from Lemma \ref{Lem:Blowup2} in this case is: $\exists$ $t_0>0$ and $R>0$ such that $f_R(t_0)>M_0=M_0(d,p,C)={}$a constant. 
Then we proceed by assuming the contrary, i.e.~$u$ is global, hence $f_R(t)\leq M_0$ for all $t>0,R>0$, and then derive a contradiction. We apply the \textit{test function method} as in \cite{MelianEtal10} (see also \cite{MitidieriEtal01}) to get the desired contradiction, thereby proving Theorem \ref{Thm:Blow-up} (2). \medskip

For Theorem \ref{Thm:Blow-up} (3), we use (\ref{Est:Bern}) to get the blow-up criterion according to Lemma \ref{Lem:Blowup2} as
\begin{align}\label{Est:bsup}
\int_{\R^n}\left(1+\frac{|x|^2}{R}\right)^{-b}u_0(x)dx>MR^{\gamma}\quad(\gamma:=n-\frac{\sigma+2}{p-1}>0).
\end{align}
Both sides of this inequality are increasing with respect to $R$ and, by the monotone convergence theorem, the left hand side converges to $\|u_0\|_{L^1}$ as $R\to\infty$. Now setting $R=2$ (in view of Theorem \ref{Thm:AsympPsi}) the right hand side of (\ref{Est:bsup}) achieves the minimum, hence Theorem \ref{Thm:Blow-up} (3) follows.

\begin{remark}\rm

Let us compare our method to the principal eigenfunction method. For $R>0$, the principal eigenvalue $\L_R$ and eigenfunction $\psi_R$ are solutions of
\[
J\ast\psi_R-\psi_R=-\lambda_R\psi_R\quad\mbox{in $B_R(x_0)$},\quad\psi_R=0\quad\mbox{on $\R^n\setminus B_R(x_0)$}.
\]
Then it was shown in \cite{MelianEtal10} that $R^2\lambda_R\to\mbox{a constant}>0$ and $R^N\psi_R(Rx)\to\Psi$ in $L^2(B_1(x_0))$, 
as $R\to\infty$, where $\Psi>0$ is the principal eigenfunction of $-\Lap$ on $B_1(x_0)$. However, this re-scaled property requires $J$ to have compact support \cite{MelianEtal09}. Taking $R$ large, we see that $\psi_R$ satisfies
\[
J\ast\psi_R-\psi_R\geq-\frac{d}{R^2}\psi_R,
\]
which is precisely the left inequality of (\ref{Ineq:LT}) with $\eta=R^2$ assuming $\A_0=\|J\|_{L^1}=1$.

\end{remark}

We end the paper with the following result.

\begin{corollary}

Assume $a(x,t)\sim\lng x\rng^\sigma$ ($\sigma\geq0$). Eq.\ (\ref{Eqn:main}) has the Fujita critical exponent
\[
p_F=1+\frac{\S+2}{n}.
\]

\end{corollary}

The rest of the paper is organized as follows. We present in Section \ref{Sec:prelim} notation, basic estimates, and some definitions. The weighted norm estimate (Theorem \ref{Thm:GreenFar}) is proved in Section \ref{Sec:WNorm} and the interpolation (Theorem \ref{Thm:Interp}) is proved in Section \ref{Sec:Interp}. In Section \ref{Sec:LocWell}, we establish the local well-posedness (Theorem \ref{Thm:LocWell}) of the Cauchy problem (\ref{Eqn:main}). The blow-up phenomena (Theorem \ref{Thm:Blow-up}) in the sub-critical, critical, and super-critical exponent are investigated in Section \ref{Sec:Blowup}. Finally, the global well-posedness is derived in Section \ref{Sec:Global}.

\section{Preliminaries}\label{Sec:prelim}

\subsection{Notation, basic results, and definitions}

\begin{itemize}
\item Weighted Lebesgue spaces
\begin{align*}
&L^p_b(\R^n)=\left\{f:\R^n\to\R\mid f\,\,\mbox{measurable},\|f\|_{L^p_b}:=\|\langle\cdot\rangle^bf\|_{L^p}<\infty\right\},\\
&L^p_\infty(\R^n):=\bigcap_{b\in\R}L^p_b(\R^n)\quad(1\leq p\leq\infty).
\end{align*}
\item Auxiliary functions
\begin{align*}
\r=\left(1+\frac{|x|^2}{\eta}\right)^{1/2}\quad\mbox{and}\quad\varGamma=\left(1+\frac{|x|^2}{\eta}\right)^{b/2}\quad(\eta>0,\,b\in\R).
\end{align*}
\end{itemize}





\begin{lemma}\label{Lem:Ineq1}

For all $a,b\in\R$, $|a^p-b^p|\leq C_p(|a|\vee|b|)^{p-1}|a-b|$.

\end{lemma}

The following result extends the well-known inequality $\lng x\rng^b\leq2^{|b|/2}\lng x-y\rng^{|b|}\lng y\rng^b$.

\begin{lemma}\label{Lem:twoside}
Let $b\in\R$, $\eta\geq2$, $r,s\in[0,1]$, and denote $\r(x)=(1+|x|^2/\eta)^{1/2}$. Then
\begin{align}\label{Aux:2side2}
\left(\frac{r\r(y)^2+(1-r)\r(x)^2}{s\r(y)^2+(1-s)\r(x)^2}\right)^{b/2}\leq2^{|b|/2}\lng x-y\rng^{|b|},
\end{align}
for all $x,y\in\R^n$.

\end{lemma}

\noindent\textbf{Proof.} 
It suffices to prove for $b>0$ that
\begin{align*}
\CR:=\left(\frac{r\r(y)^2+(1-r)\r(x)^2}{s\r(y)^2+(1-s)\r(x)^2}\right)^{b/2}\leq
\begin{cases}
\displaystyle
2^{b/2}&\mbox{if $|x-y|\leq1$},\\
\vspace{-10pt}\\
\displaystyle
\lng x-y\rng^{b}&\mbox{if $|x-y|\geq1$}.
\end{cases}
\end{align*}
Observe that if $0\leq\L\leq1$ then $|x|+|y|\leq|x-y|+2\L|y|+2(1-\L)|x|$, hence
\begin{align}
|x|+|y|\leq|x-y|+1+\L|y|^2+(1-\L)|x|^2.\label{Aux:lema2side1}
\end{align}

If $|x-y|\leq 1$, we have by (\ref{Aux:lema2side1}) that $\eta+s|y|^2+(1-s)|x|^2\geq|x|+|y|$, so
\begin{align*}
\CR&=\left(1+\frac{(r-s)(|y|^2-|x|^2)}{\eta+s|y|^2+(1-s)|x|^2}\right)^{b/2}
\leq 2^{b/2}.
\end{align*}
Assume $|x-y|\geq 1$. We have
\begin{align*}
(1+&|x-y|^2)(\eta+s|y|^2+(1-s)|x|^2)-(\eta+r|y|^2+(1-r)|x|^2)\\
&=(s-r)(|y|^2-|x|^2)+\{\eta|x-y|^2+s|x-y|^2|y|^2+(1-s)|x-y|^2|x|^2\}\\
&\geq-|s-r||x-y|(|x|+|y|)+|x-y|\{2|x-y|+s|y|^2+(1-s)|x|^2\}\\
&\geq-|x-y|(|x|+|y|)+|x-y|(|x|+|y|)\geq0.
\end{align*}
This implies $\CR\leq\lng x-y\rng^b$.\QED\medskip

We will need the sharp Young's inequality or Brascamp-Lieb inequality \cite{BrascampLieb76,LiebLoss01}.

\begin{proposition}\label{Prop:SharpYoung}
Let $k\geq2$ and $1\leq p_1,\ldots,p_k,r\leq\infty$ be such that 
\[
\frac{1}{p_1}+\cdots+\frac{1}{p_k}=k-1+\frac{1}{r}.
\]
 Then
\begin{align}\label{Est:sharpY}
\|f_1\ast\cdots\ast f_k\|_{L^r}\leq\left(\prod_{l=1}^kC_{p_l}\right)^n\|f_1\|_{L^{p_1}}\cdots\|f_k\|_{L^{p_k}}
\end{align}
for all $f_i\in L^{p_i}(\R^n)$ ($i=1,\ldots,k$), where $C_p=\sqrt{\frac{p^{1/p}}{q^{1/q}}}$ with $q=p'=p/(p-1)$.

\end{proposition}

\begin{lemma}
\label{Lem:Jkrho}

Let $\eta>0$ and $x_0,x_1,\ldots,x_k\in\R^n$. Denote $\r(x,\eta)=(1+|x|^2/\eta)^{1/2}$. Then
\begin{align}\label{Est:Tech1}
\r(x_0,\eta)^b\leq
\begin{cases}
\displaystyle
\prod_{j=1}^k\r\left(x_{j-1}-x_j,\frac{\eta}{2}\right)^b
\cdot\r\left(x_k,\frac{\eta}{2}\right)^b
&\mbox{if $b\geq0$}\\
\vspace{-8pt}\\
\displaystyle
\prod_{j=1}^k\r(x_{j-1}-x_{j},\eta)^{|b|}\cdot\r(x_k,2\eta)^b
&\mbox{if $b<0$}.
\end{cases}
\end{align}

\end{lemma}

\noindent\textbf{Proof.} 
If $b\geq0$, the desired inequality is equivalent to that
\begin{align*}
1+\frac{|x_0|^2}{\eta}&\leq\prod_{j=1}^k\left(1+\frac{|x_{j-1}-x_j|^2}{\eta/2}\right)\cdot\left(1+\frac{|x_k|^2}{\eta/2}\right)\\
&=1+\frac{2}{\eta}(|x_0-x_1|^2+\cdots+|x_{k-1}-x_k|^2+|x_k|^2)+\cdots,
\end{align*}
which is obviously true by the triangle inequality and the AM-GM inequality.\medskip

In the case $b<0$, the desired inequality is equivalent to
\[
\r(x_k,2\eta)^{-b}\leq\r(x_k-x_{k-1},\eta)^{|b|}\cdots\r(x_1-x_0,\eta)^{|b|}\r(x_0,\eta)^{-b},
\]
which is true by the previous case.
\QED\medskip

We will need the following proposition on asymptotic behavior of a regular varying coefficients exponential series, that is a series of the form
\[
\sum_{k=N}^\infty R(k)\frac{t^k}{k!},
\]
where $R$ is a regular varying function. 
  A positive measurable function $R$, defined on $(s_0,\infty)$, is called \textit{regular varying} if there is $\delta\in\R$ such that
\[
\lim_{s\to\infty}\frac{R(\L s)}{R(s)}=\L^\delta\quad(\forall\,\L>0).
\]
$\delta$ is called the index of $R$. 
If $\delta=0$, $R$ is called \textit{slowly varying}. Regular varying functions are precisely functions that can be expressed as $R(s)=s^\delta L(s)$ for some $\delta\in\R$, where $L$ is a slowly varying function.

\begin{proposition}[See \cite{BinghamEtal83,Khomrutai17}]\label{Prop:Kummer}

Let $N\in\BN$, $b\in\R$, and $L$ a slowly varying function. Then
\begin{align}\label{Est:Regv}
\sum_{k=N}^\I k^{b}L(k)\frac{t^k}{k!}\asymp t^{b}L(t)e^{t},
\end{align}
for all $t\geq t_0>0$.

\end{proposition}

\begin{lemma}\label{Lem:kernel}

$\CJ^k$ ($k\geq2 $) is an integral operator with the kernel $J_k$ given by
\[
J_k(x)=\int_{(\R^n)^{k-1}}J(x-y_1)J(y_1-y_2)\cdots J(y_{k-2}-y_{k-1}) J(y_{k-1})dy_{k-1}\cdots dy_1.
\]

\end{lemma}

\noindent\textbf{Proof.} 
The proof is
straightforward.
\QED\medskip

Using $J_k$, the Green operator $\CG(t)$ in (\ref{GreenDef}) has the kernel given by
\begin{align}\label{Exp:Gtker}
G(x,t)=e^{-\A_0t}\sum_{k=0}^\infty\frac{t^k}{k!}J_k(x).
\end{align}






\begin{definition}\label{Def:Sol}

A function $u$ is called a \textbf{(mild) solution} of (\ref{Eqn:main}) provided it satisfies
\begin{align}\label{Exp:Mild}
u(t)=\CG(t)u_0+\int_0^t\CG(t-\tau)\{a(\tau)u(\tau)^p\}d\tau
\end{align}
for all $x,t$.

\end{definition}

\subsection{$\varepsilon$-Equilibriums}\label{subsec:epequi}

The concept of equilibriums (or stationary solutions) is  important in the study of evolution equations. In particular, it leads to the so-called general relative entropy inequalities \cite{ChafaiEtal04,MichelEtal04, MichelEtal05}. 
Consider the scattering (or linear Boltzmann) equation
\begin{align}\label{Eqn:Equi}
\P_tu=\int_{\R^n}J(x,y)u(y,t)dy-m(x)u(x,t),
\end{align}
where $J$ is a given function and $m(x)=\int J(y,x)dy$.

\begin{definition}
A function $\varGamma=\varGamma(x)$ is called an \textit{equilibrium} for Eq.~(\ref{Eqn:Equi}) if it satisfies
\[
\int_{\R^n}J(x,y)\varGamma(y)dy-m(x)\varGamma(x)=0\quad(\forall\,x\in\R^n).
\]

\end{definition}

\begin{lemma}
If $\varGamma>0$ is an equilibrium for (\ref{Eqn:Equi}), then for any convex differentiable function $\Phi$ and any solution $u$ of (\ref{Eqn:Equi}) we have
\[
\frac{d}{dt}\int_{\R^n}\Phi\left(\frac{u}{\varGamma}\right)\varGamma dx\leq0.
\]

\end{lemma}

\noindent\textbf{Proof.} 
By convexity, $\Phi'(a)(a-b)-\Phi(a)+\Phi(b)\geq0$ for all $a,b\in\R$ so
\begin{align*}
&\frac{d}{dt}\int\Phi\left(\frac{u}{\varGamma}\right)\varGamma dx\\
&\,\,=-\iint J(x,y)\varGamma(y)\left\{\Phi'\left(\frac{u(x,t)}{\varGamma(x)}\right)\left(\frac{u(x,t)}{\varGamma(x)}-\frac{u(y,t)}{\varGamma(y)}\right)-\Phi\left(\frac{u(x,t)}{\varGamma(x)}\right)+\Phi\left(\frac{u(y,t)}{\varGamma(y)}\right)\right\}dydx,
\end{align*}
which implies the desired inequality.\QED\medskip

In this work, $J$ is radially symmetric, i.e.~$J=J(|x-y|)$, with $m=\|J\|_{L^1}=\A_0$, and equilibriums can be easily seen to be $\varGamma={}$a constant. However, constant functions can hardly lead to interesting results for (\ref{Eqn:main}).
So we propose the following notion.

\begin{definition}
\label{Def:Lyap}

A function
$\varGamma=\varGamma(x;\varepsilon)$ is called an \textbf{$\varepsilon$-equilibrium} for Eq.~(\ref{Eqn:Equi}) if it satisfies
\[
\left|\int_{\R^n}J(x,y)\varGamma(y;\varepsilon)dy-m(x)\varGamma(x;\varepsilon)\right|\leq\varepsilon\varGamma(x;\varepsilon)\quad(\forall\,x\in\R^n).
\]

\end{definition}


As a first application, we prove the following relative entropy inequality.

\begin{proposition}

Let $J\geq0$ be radially symmetric, $J\in L^1_\infty(\R^n)$, $m=\|J\|_{L^1}=\A_0$, and 
\[
\varGamma(x,t)=\left(1+\frac{|x|^2}{\eta_0+t}\right)^{b/2}
\]
where $\eta_0\geq2$ and $b\in\R$. Assume $\Phi\geq0$ is a convex differentiable function satisfying $\Phi(0)=0$ and $u$ a solution of (\ref{Eqn:Equi}). Then there is a constant $\nu>0$, independent of $u$, such that
\begin{align*}
\frac{d}{dt}\int_{\R^n}\Phi\left((1+t)^{-\nu}\frac{u}{\varGamma}\right)\varGamma dx\leq0.
\end{align*}

\end{proposition}

\noindent\textbf{Proof.} 
By Theorem \ref{Thm:AsympPsi}, $\varGamma$ is an $\varepsilon$-equilibrium for (\ref{Eqn:Equi}), where $\varepsilon=\frac{d}{\eta_0+t}$ ($d$ is a constant depending only upon $b$ and $J$). 
Since $\Phi$ is convex and $\Phi(0)=0$, $\Phi(0)\geq\Phi(s)+\Phi'(s)(-s)$ which implies
\[
\Phi'(s)s\geq\Phi(s)\geq0\quad\forall\,s\in\R.
\]

Let us denote $F=\L(t)u(x,t)/\varGamma(x,t)$, where $\L(t)=(1+t)^{-\nu}$. Consider
\begin{align*}
\CR&:=
\frac{d}{dt}\int\Phi(F)\varGamma dx=\int\Phi'(F)\left[\frac{\L\P_tu+\L'u}{\varGamma}-\L u\frac{\P_t\ln\varGamma}{\varGamma}\right]\varGamma dx+\int\Phi(F)\P_t\varGamma dx\\
&=\int\Phi'(F)\L\left[J\ast u-(1-(\ln\L)')u-u\P_t\ln\varGamma\right]dx+\int\Phi(F)\P_t\varGamma dx.
\end{align*}
The first term of $\CR$ equals
\begin{align*}
\int\Phi'(F)&\L J\ast u\,dx=\iint J(x-y)\Phi'(F(x,t))\varGamma(y,t)F(y,t)\,dydx\\
&=-\iint J(x-y)\varGamma(y,t)\left\{\Phi'(F(x,t))[F(x,t)-F(y,t)]\right\}dydx\\
&\hspace{4cm}+\iint J(x-y)\varGamma(y,t)\Phi'(F(x,t))F(x,t)dydx\\
&\leq-\iint J(x-y)\varGamma(y,t)\left\{\Phi'(F(x,t))[F(x,t)-F(y,t)]\right\}\,dydx\\
&\hspace{5cm}+(1+\frac{d}{\eta_0+t})\int\Phi'(F(x,t))\L(t)u(x,t)\,dx,
\end{align*}
where we have used that $\Phi'(s)s\geq0$ and $\varGamma$ is an $\varepsilon$-equilibrium.
Now we have
\begin{align*}
\CR&\leq-\iint J(x-y)\varGamma(y,t)\left\{\Phi'(F(x,t))[F(x,t)-F(y,t)]\right\}dydx\\
&\hspace{3cm}-\int\Phi'(F(x,t))\left(-\frac{d}{\eta_0+t}-(\ln\L)'+\P_t\ln\varGamma\right)\L u\,dx\\
&\hspace{5cm}+\int\Phi(F(x,t))\P_t\varGamma\, dx=:-\CR_1-\CR_2+\CR_3.
\end{align*}

Using that $\varGamma$ is $\varepsilon$-equilibrium and the Fubini's theorem, we get
\begin{align*}
&\iint J(x-y)\varGamma(y,t)\Phi(F(x,t))dydx\geq(1-\frac{d}{\eta_0+t})\int\varGamma(x,t)\Phi(F(x,t))dx\\
&\iint J(x-y)\varGamma(y,t)\Phi(F(y,t))dydx=\int\varGamma(x,t)\Phi(F(x,t))dx,
\end{align*}
which imply
\begin{align*}
&\iint J(x-y)\varGamma(y,t)\{\Phi(F(x,t))-\Phi(F(y,t))\}dydx+\frac{d}{\eta_0+t}\int\varGamma(x,t)\Phi(F(x,t))dx\\
&\hspace{3cm}=:-\CR_4+\CR_5\geq0.
\end{align*}

Combining the above estimates we obtain that
\begin{align*}
\CR\leq-(\CR_1+\CR_4)-(\CR_2-\CR_3-\CR_5).
\end{align*}
Since $\Phi$ is convex, it follows that $\Phi'(a)(a-b)-\Phi(a)+\Phi(b)\geq0$. Hence
\begin{align*}
\CR_1+\CR_4
\geq0.
\end{align*}
We calculate
\[
|\P_t\ln\varGamma|=\left|-(b/2)\frac{|x|^2}{(\eta_0+t)^2}\frac{\eta_0+t}{\eta_0+t+|x|^2}\right|\leq\frac{|b|/2}{\eta_0+t},
\]
so we have
\begin{align*}
|\CR_3|=\left|\int\varGamma \Phi(F)\P_t\ln\varGamma\,dx\right|\leq\frac{|b|/2}{\eta_0+t}\int\varGamma\Phi(F)\,dx.
\end{align*}
Using that $\Phi'(s)s\geq\Phi(s)$ and by taking $\L=(1+t)^{-(2d+|b|)}$, we get
\begin{align*}
\CR_2-\CR_3-\CR_5&\geq\int\Phi'(F)\left(-\frac{d}{\eta_0+t}-(\ln\L)'+\P_t\ln\varGamma\right)\L u\,dx-\frac{d+|b|/2}{\eta_0+t}\int\varGamma\Phi(F)\,dx\\
&\geq\int\Phi'(F)\frac{d+|b|/2}{\eta_0+t}\L u\,dx-\frac{d+|b|/2}{\eta_0+t}\int\varGamma\Phi(F)\,dx\\
&\geq\frac{d+|b|/2}{\eta_0+t}\int\varGamma\Phi'(F)F\,dx-\frac{d+|b|/2}{\eta_0+t}\int\varGamma\Phi(F)\,dx\geq0.
\end{align*}
So the desired inequality $\CR\leq0$ is true.
\QED

\section{Proof of Theorem \ref{Thm:GreenFar}
}\label{Sec:WNorm}

\begin{lemma}\label{Lem:GBr}

Let $\varGamma=(1+|x|^2/\eta)^{b/2}$ where $\eta\geq1,x\in\R^n$. Then
\[
\eta^{-b_+/2}\lng x\rng^b\leq\varGamma\leq\eta^{-b_-/2}\lng x\rng^b,
\]
where $b_+=\max\{b,0\}$ and $b_-=\min\{b,0\}$.

\end{lemma}

\noindent\textbf{Proof.} This follows from that $b=b_++b_-$, hence $\varGamma=\eta^{-b_+/2}\eta^{-b_-/2}(\eta+|x|^2)^{b/2}$.\QED

\subsection{An $\varepsilon$-equilibrium result}

\begin{theorem}\label{Thm:AsympPsi}

Let $J\in L^{1}_{\delta_0}(\R^n)$ ($\delta_0\geq2$), $|b|\leq\delta_0-2$, and $\r(x)=(1+|x|^2/\eta)^{1/2}$ ($\eta\geq2$). Then there is a constant $d=d(b,J)>0$ such that
\begin{align}\label{Est:Lyap}
\left|\int_{\R^n}|J(x-y)|\r(y)^bdy-\A_0\r(x)^b\right|\leq\frac{d}{\eta}\r(x)^b\quad(\forall\,x\in\R^n).
\end{align}
In fact, $d=C_b\|J\|_{L^1_{2+|b|}}$ and we can take $C_b=0$ if $b=0$.

\end{theorem}

\noindent\textbf{Proof.} 
By Taylor's theorem, we have
\begin{align*}
\frac{\r(y)^b}{\r(x)^b}-1&
=\left(1+\frac{|y|^2-|x|^2}{\eta+|x|^2}\right)^{b/2}-1\\
&=\kappa S+\kappa(\kappa-1)\int_0^1S^2(1+\L S)^{\kappa-2}(1-\L)d\L,
\end{align*}
where
$S=(|y|^2-|x|^2)/(\eta+|x|^2)$ and $\kappa=b/2$. Then by integration we get
\begin{align*}
\CR&:=\int|J(x-y)|\left(\frac{\r(y)^b}{\r(x)^b}-1\right)dy\\
&=\kappa\int|J(x-y)|Sdy+\kappa(\kappa-1)\int\int_0^1|J(x-y)|S^2(1+\L S)^{\kappa-2}(1-\L)d\L dy\\
&=:\CR_1+\CR_2.
\end{align*}

To estimate $\CR_1$, we note that $\int_{\R^n}|J(z)|zdz=0$. Then 
\begin{align*}
&S=\frac{|x-y|^2}{\eta+|x|^2}-\frac{2x\cdot(x-y)}{\eta+|x|^2}\quad\Rightarrow\quad
|\CR_1|
\leq\left(\frac{|b|}{2}\|J\|_{L^1_2}\right)\frac{1}{\eta}.
\end{align*}
Next we estimate $\CR_2$. Note that\begin{align*}
S^2(1+\L S)^{\kappa-2}
&=\left(\frac{\L\r(y)^2+(1-\L)\r(x)^2}{\r(x)^2}\right)^{b/2}\left(\frac{|y|^2-|x|^2}{\eta+\L|y|^2+(1-\L)|x|^2}\right)^2.
\end{align*}
By Lemma \ref{Lem:twoside}, we have
\begin{align*}
\left(\frac{\L\r(y)^2+(1-\L)\r(x)^2}{\r(x)^2}\right)^{b/2}&\leq2^{|b|/2}\lng x-y\rng^{|b|}
\end{align*}
then by direct integration we get
\begin{align*}
\int_0^1|J(x-y)|S^2(1+\L S)^{\kappa-2}(1-\L)d\L
=\frac{(|y|^2-|x|^2)^2}{(\eta+|x|^2)(\eta+|y|^2)}.
\end{align*}
Thus
\begin{align*}
|\CR_2|&\leq\frac{|b(b-2)|}{4}2^{|b|/2}\int|J(x-y)|\lng x-y\rng^{|b|}\frac{(|y|^2-|x|^2)^2}{(\eta+|x|^2)(\eta+|y|^2)}dy.
\end{align*}
Expanding
\begin{align*}
(|y|^2-|x|^2)^2&=|x-y|^4-4(x-y)\cdot x+4((x-y)\cdot x)^2,
\end{align*}
and using $|x-y|^2\leq\frac{2}{\eta}(\eta+|x|^2)(\eta+|y|^2)$ and $|x|\leq\frac{1}{2\eta^{3/2}}(\eta+|x|^2)(\eta+|y|^2)$, we get
\begin{align*}
|\CR_2|&\leq C_b\|J\|_{L^1_{2+|b|}}\frac{1}{\eta}.
\end{align*}

Combining the estimates of $\CR_1,\CR_2$ above, we obtain
\[
\left|\int|J(x-y)|\left(\frac{\r(y)^b}{\r(x)^b}-1\right)dy\right|\leq \frac{d}{\eta},
\]
where $d=C_b\|J\|_{L^1_{2+|b|}}$. This implies the result.\QED\medskip

\begin{remark}\rm

If $J$ has algebraic tail
\[
J(x)\sim\frac{1}{|x|^{n+r}}\quad(0<r<2),\,\mbox{as $|x|\to\infty$},
\]
then it does not satisfy the hypothesis of Theorem \ref{Thm:AsympPsi}. 
This kernel occurs in the study of stable laws with index $r$. For interesting recent results on the nonlocal diffusion equation in this case, especially the Fujita critical exponent, see \cite{Alfaro17}.
\end{remark}

\subsection{Proof of Theorem \ref{Thm:GreenFar}}

Let $f\in L^q_b(\R^n)$. We split the proof into the cases $q=1$ and $q=\infty$.\smallskip

\textbf{Case $q=1$.} Let $\eta_0\geq2$. We denote 
\[
\varGamma(x,t)=\left(1+\frac{|x|^2}{\eta_0+t}\right)^{b/2},\quad\eta_t=\eta_0+t.
\] 
Below, $\R^n$ will be omitted in writing all integrals over this set, and similar for $(\R^n)^k$. By Theorem \ref{Thm:AsympPsi}, we get
\begin{align}
\int|J(x-y)|\varGamma(y,t)dy\leq\left(\A_0+\frac{d}{\eta_t}\right)\varGamma(x,t).\label{Temp1:ThmGreenFar}
\end{align}
Using the triangle inequality, Lemma \ref{Lem:GBr}, the Fubini's theorem, and (\ref{Temp1:ThmGreenFar}), we have
\begin{align*}
\|\CJ^kf\|_{L^1_b}&\leq\int\lng x\rng^b\int|J(x-y_1)\cdots J(y_{k-1}-y_k)f(y_k)|dy_k\cdots dy_1dx\\
&\leq\eta_t^{b_+/2}\int\int\varGamma(x,t)|J(x-y_1)|dxJ(y_1-y_2)\cdots J(y_{k-1}-y_k)|f(y_k)|dy_k\cdots dy_1\\
&\leq\eta_t^{b_+/2}\left(\A_0+\frac{d}{\eta_t}\right)\int\varGamma(y_1,t)|J(y_1-y_2)\cdots J(y_{k-1}-y_k)f(y_k)|dy_k\cdots dy_1.
\end{align*}
Applying (\ref{Temp1:ThmGreenFar}) repeatedly to integrate $dy_1,dy_2,\ldots,dy_{k-1}$, respectively, we then get
\begin{align*}
\|\CJ^kf\|_{L^1_b}&\leq\eta_t^{b_+/2}\left(\A_0+\frac{d}{\eta_t}\right)^k\int\varGamma(y_k,t)|f(y_k)|dy_k.
\end{align*}
Using Lemma \ref{Lem:GBr} and recalling $\eta_t=\eta_0+t$, we obtain
\begin{align*}
\|\CJ^kf\|_{L^1_b}
&\leq(\eta_0+t)^{|b|/2}\left(\A_0+\frac{d}{\eta_0+t}\right)^k\|f\|_{L^1_b}.
\end{align*}

The above estimate is true for all $k\geq1$ and so is $k=0$, hence
\begin{align*}
\|\CG(t)f\|_{L^1_b}&\leq e^{-\A_0t}\sum_{k=0}^\I\frac{t^k}{k!}\|\CJ^kf\|_{L^1_b}\\
&\leq e^{-\A_0t}(\eta_0+t)^{|b|/2}\|f\|_{L^1_b}\sum_{k=0}^\I\frac{t^k}{k!}\left(\A_0+\frac{d}{\eta_0+t}\right)^k\\
&\leq(\eta_0+t)^{|b|/2}\|f\|_{L^1_b}\exp\left(-\A_0t+t\left(\A_0+\frac{d}{\eta_0+t}\right)\right).
\end{align*}
It follows that 
\begin{align*}
\|\CG(t)f\|_{L^1_b}
&\lesssim \lng t\rng^{|b|/2}\|f\|_{L^1_b}.
\end{align*}
This proves (\ref{Est:ThmGreenFar1}) for the case $q=1$.\medskip

\textbf{Case $q=\infty$.} We shall prove the estimate (\ref{Est:ThmGreenFar1}) in the case $q=\I$. We put
\[
\varGamma(x,t)=\left(1+\frac{|x|^2}{\eta_t}\right)^{-b/2}.
\]
Substituting $-b$ for $b$ in Lemma \ref{Lem:GBr} and using $(-b)_+=-b_-,(-b)_-=-b_+$, we get
\[
\eta_t^{b_-/2}\lng x\rng^{-b}\leq\varGamma(x,t)\leq\eta_t^{b_+/2}\lng x\rng^{-b}.
\]
For each $x\in\R^n$, we have
\begin{align*}
|\CJ^kf(x)|&=\left|\int J(x-y_1)\cdots J(y_{k-1}-y_k)\lng y_k\rng^{-b}\lng y_k\rng^bf(y_k)dy_k\cdots dy_1\right|\\
&\leq\|f\|_{L^\infty_b}\int|J(x-y_1)\cdots J(y_{k-1}-y_k)|\lng y_k\rng^{-b}dy_k\cdots dy_1\\
&\leq\eta_t^{-b_-/2}\|f\|_{L^\infty_b}\int|J(x-y_1)\cdots J(y_{k-1}-y_k)|\varGamma(y_k,t) dy_k\cdots dy_1,\\
&\leq\eta_t^{-b_-/2}\|f\|_{L^\infty_b}\left(\A_0+\frac{d}{\eta_t}\right)\int|J(x-y_1)\cdots J(y_{k-2}-y_{k-1})|\varGamma(y_{k-1},t)dy_{k-1}\cdots dy_1,
\end{align*}
where we have used (\ref{Temp1:ThmGreenFar}). By repeatedly integrating in $dy_1,\ldots,dy_{k-1}$ and employing the same argument, we get
\begin{align*}
|\CJ^kf(x)|
&\leq\eta_t^{-b_-/2}\|f\|_{L^\infty_b}\left(\A_0+\frac{d}{\eta_t}\right)^k\varGamma(x,t).
\end{align*}
Now recalling that $\varGamma(x,t)\leq\eta^{b_+/2}\lng x\rng^{-b}$, hence we obtain
\begin{align*}
\|\CJ^kf\|_{L^\infty_b}=\sup_x|\lng x\rng^b\CJ^kf(x)|&\leq\eta_t^{|b|/2}\|f\|_{L^\infty_b}\left(\A_0+\frac{d}{\eta_t}\right)^k.
\end{align*}
We estimate the Green operator
\begin{align*}
\|\CG(t)f\|_{L^\I_b}&\leq e^{-\A_0t}\sum_{k=0}^\I\frac{t^k}{k!}\|\CJ^kf\|_{L^\infty_b}\\
&\leq e^{-\A_0t}(\eta_0+t)^{|b|/2}\|f\|_{L^\infty_b}\sum_{k=0}^\I\frac{t^k}{k!}\left(\A_0+\frac{d}{\eta_0+t}\right)^k\\
&\lesssim \lng t\rng^{|b|/2}\|f\|_{L^\infty_b},
\end{align*}
hence (\ref{Est:ThmGreenFar1}) is proved.\medskip

Finally, applying the Riesz-Thorin interpolation theorem we conclude that
\[
\|\CG(t)f\|_{L^q_b}\lesssim \lng t\rng^{|b|/2}\|f\|_{L^q_b}\quad(t\geq0).\QED
\]

\section{Proof of Theorem \ref{Thm:LocWell}}\label{Sec:LocWell}

We fix $b\in\R$, $T>0$ to be chosen, and denote the Banach spaces
\begin{align*}
&\CZ:=L^\infty_b(\R^n)\cap C(\R^n)\quad\mbox{with the norm}\\
&\qquad\|f\|:=\|f\|_{L^\infty_b},\\
&\CZ_T:=C\left([0,T];\CZ\right)\quad\mbox{with the norm}\\
&\qquad\|v\|_{\CZ_T}=\sup_{0\leq t\leq T}\|v(t)\|=\sup_{x\in\R^n,0\leq t\leq T}\lng x\rng^b|v(x,t)|.
\end{align*}
According to Definition \ref{Def:Sol} a solution of (\ref{Eqn:main}) is a fixed point of the operator
\begin{align*}
\CM v(t):=\CG(t)u_0+\int_0^t\CG(t-\tau)\{a(\tau)v(\tau)^p\}d\tau\quad(0\leq t\leq T).
\end{align*}

First we show that $\CM:\CZ_T\to\CZ_T$. By (\ref{Hyp:LocH1}) and (\ref{Hyp:LocH3}), we have
\begin{align}\label{WTrick}
\lng x\rng^b|a(x,t)|\lesssim\lng x\rng^{b+\sigma}\lesssim\lng x\rng^{bp}.
\end{align}
Also, $|b|\leq\delta-2$ (see (\ref{Hyp:LocH3})), it follows by Theorem \ref{Thm:GreenFar} 
that, $\CG(t):\CZ\to\CZ$ is a bounded linear operator and $\|\CG(t)f\|\leq\lng t\rng^{|b|/2}\|f\|$. 
Let $v\in\CZ_T$. By the triangle inequality, we estimate
\begin{align*}
\|\CM v(t)\|&\leq\left\|\CG(t)u_0\right\|+\int_0^t\|\CG(t-\tau)\{a(\tau)v(\tau)^p\}\|d\tau\\
&\leq C\lng t\rng^{|b|/2}\|u_0\|+C\int_0^t\lng t-\tau\rng^{|b|/2}\|a(\tau)v(\tau)^p\|d\tau\\
&\leq C\lng t\rng^{|b|/2}\|u_0\|+C\lng t\rng^{|b|/2}\int_0^t\sup_{x,\tau}|\lng x\rng^{b}a(x,\tau)v(x,\tau)^p|d\tau\\
&\leq C\lng t\rng^{|b|/2}\|u_0\|+C\lng t\rng^{|b|/2}\int_0^t\sup_{x,\tau}\lng x\rng^{bp}|v(x,\tau)|^pd\tau\qquad(\mbox{by (\ref{WTrick})})\\
&\leq C\lng t\rng^{|b|/2}\|u_0\|+C\lng t\rng^{|b|/2}t\|v\|_{\CZ_T}^p.
\end{align*}
Thus $\CM v(t)\in L^\infty_b(\R^n)$ for all $0\leq t\leq T$. It is clearly that $\CM v(t)\in C(\R^n)$ so $\CM v(t)\in\CZ$. Also, by the semigroup property of $\CG(t)$, we have
\begin{align*}
\|\CM v(s)-\CM v(t)\|&\leq\|(\CG(s)-\CG(t))u_0\|+\int_s^t\|\CG(t-\tau)\{a(\tau)v(\tau)^p\}\|d\tau\\
&\hspace{2.8cm}+\int_0^s\|(\CG(t-\tau)-\CG(s-\tau))\{a(\tau)v(\tau)^p\}\|d\tau,
\end{align*}
which implies $\CM v(s)\to\CM v(t)$ in $\CZ$, as $s\to t$. Hence
\begin{align}\label{Est:LocTmp1}
\begin{cases}
\CM:\CZ_T\to\CZ_T\quad\mbox{and}\\
\vspace{-10pt}\\
\|\CM v\|_{\CZ_T}\leq C\lng T\rng^{|b|/2}\left(\|u_0\|+T\|v\|_{\CZ_T}^p\right). 
\end{cases}
\end{align}

Next we introduce
\[
K>0\quad\mbox{and}\quad\CB_{T}(K):=\{v\in\CZ_T:\|v\|_{\CZ_T}\leq K\}.
\]
$\CB_T(K)$ is a complete metric space. We are going to show that $\CM$ is a self-map on $\CB_T(K)$ for $K$ large enough and $T$ small enough. By (\ref{Est:LocTmp1}) and $\|v\|_{\CZ_T}\leq K$, we have
\begin{align*}
\|\CM v\|_{\CZ_T}\leq C\lng T\rng^{|b|/2}\left(\|u_0\|+TK^p\right).
\end{align*}
Choose a pair $K,T$ satisfying
\begin{align*}
\lng T\rng^{|b|/2}\|u_0\|\leq \frac{K}{2C},\quad C\lng T\rng^{|b|/2}T\leq\frac{1}{2K^{p-1}}.
\end{align*}
For $v,w\in\CB_T(K)$, we have
\begin{align*}
\|\CM v(t)-\CM w\|_{\CZ_T}&\leq\sup_{x,t}\int_0^t\left\|\CG(t-\tau)\left\{a(\tau)\left(v(\tau)^p-w(\tau)^p\right)\right\}\right\|d\tau\\
&\leq C\lng T\rng^{|b|/2}\sup_{x,t}\int_0^t\left\|a(\tau)\left(v(\tau)^p-w(\tau)^p\right)\right\|d\tau\\
&\leq C\lng T\rng^{|b|/2}TK^{p-1}\|v-w\|_{\CZ_T}.
\end{align*}
Taking $T$ even smaller, it follows that $\CM:\CB_T(K)\to\CB_T(K)$ is a contraction. By the Banach fixed point theorem, we conclude that $\CM$ has a unique fixed point in $\CB_T(K)$, and hence (\ref{Eqn:main}) has a unique solution in $\CZ_T$, for $T>0$ small. This proves the first part of Theorem \ref{Thm:LocWell}.\medskip

Next, we prove the positivity of solutions assuming $J\geq0,a\geq0$. In this case we consider the positive cone 
\[
\CB_T(K)^+:=\{v\in\CB_T(K):v\geq0\}.
\]
If $v\in\CB_T(K)^+$ then $\CM v\geq0$, hence $\CM$ is a self-map on $\CB_T(K)^+$ for $K$ large and $T>0$ small. By the same argument as above, it follows that there is a unique fixed point of $\CM$ in $\CB_T(K)^+$. Therefore there is a unique positive solution $u$.\QED

\section{Proof of Theorem \ref{Thm:Interp}}\label{Sec:Interp}

Take $N$ to be a sufficiently large positive integer to be specified and decompose the series (\ref{Exp:Gtker}) into
\begin{align}\label{Gt:Decomp1}
G(x,t)=G_N(x,t)+R_N(x,t),
\end{align}
where $G_N$ is the sum of first $N$ terms of the series and the remainder series $R_N$ is
\begin{align}
R_N=e^{-\A_0t}\sum_{k=N}^\I\frac{t^k}{k!}J_k(x).\label{RN}
\end{align}
Let us denote the corresponding operators by $\CG(t)=\CG_N(t)+\CR_N(t)$.


\subsection{Pointwise estimate for $R_N$}

\begin{proposition}\label{Prop:RN}

Assume $\varepsilon_0>0$, $\beta\geq0$, and
\begin{align}\label{Assump:Jglobal}
J\in L^1_{2+\beta}(\R^n)\cap L^{1+\V_0}_{\beta}(\R^n).
\end{align}
Let $N\geq\lceil\frac{1}{\V_0}\rceil+1$. Then the following pointwise estimate holds
\begin{align}\label{Est:RN}
\left|R_N(x,t)\right|\lesssim \left\lng\frac{\lng x\rng^2}{\lng t\rng}\right\rng^{-\beta/2}\lng t\rng^{-\frac{n}{2}}\quad(\forall\,x\in\R^n,t\geq0).
\end{align}
In particular, if $J\in L^1_\infty(\R^n)\cap L^{1+\varepsilon_0}_\infty(\R^n)$ then (\ref{Est:RN}) is true for all $\beta\geq0$.

\end{proposition}

\noindent\textbf{Proof.} 
We shall employ the sharp Young's inequality (Proposition \ref{Prop:SharpYoung}). Let $0\leq b\leq\beta$ and $\eta\geq4$. We denote
\[
\varGamma:=\left(1+\frac{|x|^2}{\eta}\right)^{b/2},\quad\T\varGamma:=\left(1+\frac{|x|^2}{\eta/2}\right)^{b/2}.
\]
By Lemma \ref{Lem:Jkrho}, we have 
\[
\varGamma(x)\leq\T\varGamma(x-y_1)\T\varGamma(y_1-y_2)\cdots\T\varGamma(y_{k-2}-y_{k-1})\T\varGamma(y_{k-1}),
\]
for all $x,y_1,\ldots,y_{k-1}\in\R^n$. Then plugging into the identity in Lemma \ref{Lem:kernel}
\begin{align}
&|\varGamma(x)J_k(x)|\leq\int\left\{\prod_{j=1}^{k-1}\T\varGamma(y_{j-1}-y_j)J(y_{j-1}-y_j)\right\}\T\varGamma(y_{k-1})J(y_{k-1})dy_{k-1}\cdots dy_1,\label{Intp:tmp1}
\end{align}
where $y_0:=x$.\medskip 

Let $k\geq N_0=\lceil\frac{1}{\varepsilon_0}\rceil+1$. We shall apply Proposition \ref{Prop:SharpYoung} with the following parameters and functions:
\begin{align*}
&p_1=\cdots=p_k=\frac{k}{k-1},\quad r=\I,\\
&f_1=\cdots=f_k=\T\varGamma J.
\end{align*}
Note that $f_l\in L^{p_l}(\R^n)$ for all $l$. For each $p_l$, the H\"older conjugate is $q_l=k$. We calculate
\begin{align*}
C_{p_l}^k&
=\frac{1}{k^{1/2}}\left(1+\frac{1}{k-1}\right)^{(k-1)/2}\leq\frac{\sqrt{e}}{k^{1/2}}.
\end{align*}
Now by (\ref{Intp:tmp1}) and the sharp Young's inequality, we get
\begin{align}
\|\varGamma J_k\|_{L^\I}&\leq\|f_1\ast\cdots \ast f_k\|_{L^\I}\lesssim k^{-n/2}\left\|\T\varGamma  J\right\|^k_{L^{k/(k-1)}},\nonumber\\
\therefore\,\,\|\varGamma J_k\|_{L^\I}&\lesssim k^{-n/2}\|\T\varGamma J\|_{L^1}^k\left(\int|\T{J}(x)|^{k/(k-1)}dx\right)^{k-1},\label{Est:RNtmp1}
\end{align}
where $\T{J}=\T\varGamma J/\|\T\varGamma J\|_{L^1}$. 
Note that the right hand side of (\ref{Est:RNtmp1}) is finite because $k/(k-1)\leq1+\V_0$ for all $k\geq N_0$ and $J\in L^{1+\V_0}_\beta(\R^n)$.\medskip

Setting $\V=1/(k-1)$ and applying the L'Hopital rule and the dominated convergence theorem, we get
\begin{align*}
\ln\left(\int|\T{J}(x)|^{k/(k-1)}dx\right)^{k-1}&=\frac{1}{\V}\ln\int|\T{J}(x)|^{\V+1}dx\to\int|\T{J}(x)|\ln|\T{J}(x)|dx,
\end{align*}
as $k\to\infty$, where we have used that $\|\T{J}\|_{L^1}=1$. Now fix $N\geq N_0$. Then there is a constant $\gamma=\gamma(N,J)>0$ such that
\begin{align*}
\left(\int|\T{J}(x)|^{k/(k-1)}dx\right)^{k-1}\leq\exp\left(\gamma\int|\T{J}(x)|\ln|\T{J}(x)|dx\right),
\end{align*}
for all $k\geq N$.\medskip

Since $b\geq0$ and $\eta\geq4$, we get $|\T\varGamma J|\leq\lng x\rng^{b}|J|$ which implies $|\ln|\T{J}||\leq |b\ln\langle x\rangle+\ln|J|+\ln\|J\|_{L^1_b}|$. Hence
\[
\exp\left(\gamma\int|\T{J}(x)|\ln|\T{J}(x)|dx\right)\leq C<\infty,
\]
and the constant $C$ is independent of $\eta$. So we obtain by (\ref{Est:RNtmp1}) that
\[
\|\varGamma J_k\|_{L^\infty}\leq C\|\T\varGamma J\|_{L^1}^kk^{-n/2}.
\]
We estimate $\|\T\varGamma J\|_{L^1}$. By Theorem \ref{Thm:AsympPsi}, we get
\begin{align}
\|\T\varGamma J\|_{L^1}&=\int|J(y)|\left(1+\frac{|y|^2}{\eta/2}\right)^{\beta/2} dy
\nonumber\\
&\leq\A_0+\frac{d}{\eta/2},\nonumber\\
\therefore\,\,\,\,\|\varGamma J_k\|_{L^\I}&\leq C\left(\A_0+\frac{2d}{\eta}\right)^kk^{-n/2}\quad(\forall\,k\geq N),\label{Aux:Thm}
\end{align}
where $\A_0=\|J\|_{L^1}$.\medskip

We derive point-wise bounds for $R_N$ by setting $\eta=\eta_0+t$, so $\varGamma=(1+\frac{|x|^2}{\eta_0+t})^{b/2}$ in (\ref{Aux:Thm}), with two values for the exponent $b$.\medskip

\textbf{Bound I.} Using (\ref{Aux:Thm}) with $b=0$, i.e.~$\varGamma=1$, we have that
\begin{align*}
\|J_k\|_{L^\infty}\leq C_1\A_0^kk^{-n/2}\quad(\forall\,k\geq N),
\end{align*}
so by Proposition \ref{Prop:Kummer} we get
\begin{align}\label{Est:RN1}
|R_N(x,t)|&\leq C_1e^{-\A_0t}\sum_{k=N}^\infty\frac{(\A_0t)^k}{k!}k^{-n/2}\leq C_1\lng t\rng^{-n/2},
\end{align}
for all $t>0$. 
\medskip

\textbf{Bound II.} Using (\ref{Aux:Thm}) with $b=\beta$ and applying Lemma \ref{Lem:GBr} we get
\begin{align*}
|\lng x\rng^{\beta}&R_N(x,t)|\leq(\eta_0+t)^{\beta/2}\|\varGamma R_N\|_{L^\I}\\
&\leq(\eta_0+t)^{\beta/2}e^{-\A_0t}\sum_{k=N}^\I\frac{t^k}{k!}\|\varGamma J_k\|_{L^\I}\\
&\lesssim(\eta_0+t)^{\beta/2}e^{-\A_0t}\sum_{k=N}^\I\frac{t^k}{k!}\left(\A_0+\frac{2d}{\eta_0+t}\right)^kk^{-n/2}\\
&\lesssim(\eta_0+t)^{\beta/2}e^{-\A_0t}\left\lng\left(\A_0+\frac{2d}{\eta_0+t}\right)t\right\rng^{-n/2}\exp\left(\left(\A_0+\frac{2d}{\eta_0+t}\right)t\right)\\
&\lesssim\lng t\rng^{-\frac{n}{2}+\frac{\beta}{2}}.
\end{align*}
So we have
\begin{align}
|R_N(x,t)|&\lesssim \left(\frac{\lng x\rng^2}{\lng t\rng}\right)^{-\beta/2}\lng t\rng^{-n/2}.\label{Est:RN2}
\end{align}

Combining (\ref{Est:RN1}) and (\ref{Est:RN2}) we conclude that
\[
|R_N(x,t)|\lesssim\left\lng\frac{\lng x\rng^2}{\lng t\rng}\right\rng^{-b/2}\lng t\rng^{-n/2},
\]
which is the desired estimate.\QED




\begin{corollary}

If $J\in L^1_{\infty}(\R^n)\cap L^{\infty}_\infty(\R^n)$,  
then for any $\beta\geq0$ the following estimate holds
\[
|G(x,t)-e^{-\A_0t}|\lesssim_\beta\left\lng\frac{\lng x\rng^2}{\lng t\rng}\right\rng^{-\beta/2}\lng t\rng^{-n/2}.
\]

\end{corollary}

\noindent\textbf{Proof.} 
By the assumption on $J$, we can take $N=2$ in Proposition \ref{Prop:RN} so that
\[
|R_2(x,t)|\lesssim\left\lng\frac{\lng x\rng^2}{\lng t\rng}\right\rng^{-\beta/2}\lng t\rng^{-n/2}.
\]
By the assumption $J\in L^\infty_\infty(\R^n)$ and the triangle inequality, we have
\begin{align*}
|G(x,t)-e^{-\A_0t}|&\leq e^{-\A_0t}t|J(x)|+|R_2(x,t)|\\
&\lesssim\lng x\rng^{-\beta}e^{-\A_0t}t+\left\lng\frac{\lng x\rng^2}{\lng t\rng}\right\rng^{-\beta/2}\lng t\rng^{-n/2}\\
&\lesssim\left\lng\frac{\lng x\rng^2}{\lng t\rng}\right\rng^{-\beta/2}\lng t\rng^{-n/2}.\QED
\end{align*}

\subsection{Proof Theorem \ref{Thm:Interp}}

Let $N\geq N_0:=\lceil\frac{1}{\varepsilon_0}\rceil+1$ and $R_N$ be as in (\ref{RN}). \medskip

\noindent\textbf{Claim.} 
Let $1\leq r\leq\infty$ and $b<\beta-\frac{n}{r}$. Then
\[
\|R_N(\cdot,t)\|_{L^r_{b}}\leq C\lng t\rng^{-\frac{n}{2}+\left(\frac{b}{2}+\frac{n}{2r}\right)_+}\quad\forall\,t\geq0.
\]

\noindent\textbf{Proof of Claim.} 
By (\ref{Est:InterpHypJ}), it follows from Proposition \ref{Prop:RN} that
\[
|R_N(x,t)|\leq C\lng t\rng^{-n/2}\left\lng\frac{\lng x\rng^2}{\lng t\rng}\right\rng^{-\beta/2}\quad\forall\,x\in\R^n,t\geq0.
\]
We estimate
\begin{align*}
\|R_N(\cdot,t)\|_{L^r_b}&\leq C\lng t\rng^{-n/2}\left\|\lng x\rng^b\left\lng\frac{\lng x\rng^2}{\lng t\rng}\right\rng^{-\beta/2}\right\|_{L^r}\\
&\leq C\lng t\rng^{-n/2}\left(1+\left\||x|^b\left\lng\frac{|x|^2}{\lng t\rng}\right\rng^{-\beta/2}\right\|_{L^r(|x|\geq1)}\right)\\
&\leq C\lng t\rng^{-n/2}\left(1+\left\||x|^b\right\|_{L^r\left(1\leq|x|<\lng t\rng^{1/2}\right)}+\left\||x|^b\left(\frac{|x|^2}{\lng t\rng}\right)^{-\beta/2}\right\|_{L^r\left(|x|\geq\lng t\rng^{1/2}\right)}\right)\\
&\leq C\lng t\rng^{-n/2}\left(1+\lng t\rng^{\frac{b}{2}+\frac{n}{2r}}\left\||\upsilon|^b\right\|_{L^r(|\upsilon|<1)}+\lng t\rng^{\frac{b}{2}+\frac{n}{2r}}\left\||\upsilon|^{b-\beta}\right\|_{L^r(|\upsilon|\geq1)}\right)
\end{align*}
where we have substituted $\upsilon=x\lng t\rng^{-1/2}$ in the above calculations. Since $b<\beta-\frac{n}{r}$, we have $
r(b-\beta)<-n$. So $\left\||\upsilon|^{b-\beta}\right\|_{L^r(|\upsilon|\geq1)}<\infty$, and moreover, the following estimate is true
\[
\|R_N(\cdot,t)\|_{L^r_b}\leq C\lng t\rng^{-\frac{n}{2}+\left(\frac{b}{2}+\frac{n}{2r}\right)_+}\quad\forall\,t\geq0.\QED
\]

Let $1\leq q\leq Q\leq\infty$ and $1\leq r\leq\infty$ be such that
\[
\frac{1}{Q}+1=\frac{1}{q}+\frac{1}{r}.
\]
Recall $\CR_N(t)f:=R_N(\cdot,t)\ast f$. Using the inequality $\lng x\rng^b\leq C\lng x-y\rng^{|b|}+C\lng y\rng^b$, for all $b\in\R$, we get
\begin{align*}
\|\CR_N(t)f\|_{L^{Q}_b}&=\left\|\lng x\rng^b\int_{\R^n}R_N(x-y,t)f(y)dy\right\|_{L^{Q}(dx)}\\
&\leq C\left\|\int_{\R^n}\lng x-y\rng^{|b|}R_N(x-y,t)f(y)dy+\int_{\R^n}R_N(x-y,t)\lng y\rng^bf(y)dy\right\|_{L^{\T q}(dx)},
\end{align*}
hence by Young's inequality we find that
\begin{align*}
\|\CR_N(t)f\|_{L^{Q}_b}&\leq C\|\lng \cdot\rng^{|b|}R_N(\cdot,t)\|_{L^{r}}\|f\|_{L^q}+C\|R_N(\cdot,t)\|_{L^r}\|\lng\cdot\rng^bf\|_{L^q}\\
&\lesssim\lng t\rng^{\frac{|b|}{2}+\frac{n}{2r}-\frac{n}{2}}\|f\|_{L^q}+\lng t\rng^{\frac{n}{2r}-\frac{n}{2}}\|f\|_{L^q_b}.
\end{align*}
Now $\frac{|b|}{2}+\frac{n}{2r}-\frac{n}{2}
=\frac{n}{2}\left(\frac{1}{Q}-\frac{1}{q}+\frac{|b|}{n}\right)$, 
therefore we obtain
\begin{align}\label{Interp:RN}
\|\CR_N(t)f\|_{L^{Q}_b}\lesssim\lng t\rng^{\frac{n}{2}\left(\frac{1}{Q}-\frac{1}{q}+\frac{|b|}{n}\right)}\|f\|_{L^q}+\lng t\rng^{\frac{n}{2}\left(\frac{1}{Q}-\frac{1}{q}\right)}\|f\|_{L^q_b}.
\end{align}

It remains to find the bound of $\CG_N(t)f:=(G(\cdot,t)-R_N(\cdot,t))\ast f$. By (\ref{Gt:Decomp1}) and the triangle inequality,
\begin{align*}
\|\CG_N(t)f\|_{L^Q_b}\leq e^{-t}\sum_{k=0}^{N-1}\frac{t^k}{k!}\|\CJ^kf\|_{L^{Q}_b}.
\end{align*}
To estimate the term $\CJ^kf$ we use that $\lng x\rng^b\leq C\lng x-y\rng^{|b|}\lng y\rng^b$ and apply Young's inequality to get
\begin{align*}
\|\CJ^kf\|_{L^{Q}_b}&=\left\|\lng x\rng^b\int_{\R^n}J_k(x-y)f(y)dy\right\|_{L^{Q}(x)}\\
&\leq C\|\lng\cdot\rng^{|b|}J_k\|_{L^1}\|\lng\cdot\rng^bf\|_{L^{Q}}.
\end{align*}
Since $J\in L^1_{2+\beta}(\R^n)$ and $|b|<\beta$, we have $\|\lng\cdot\rng^{|b|}J\|_{L^1}<\infty$. So $\|\lng\cdot\rng^bJ_k\|_{L^1}<\infty$ for $k=1,\ldots,N-1$. Hence
\begin{align}\label{Interp:GN}
\|\CG_N(t)f\|_{L^Q_b}\leq Ce^{-t}\sum_{k=0}^{N-1}\frac{t^k}{k!}\|f\|_{L^{Q}_b}\leq Ce^{-t/2}\|f\|_{L^{Q}_b}.
\end{align}
Combining (\ref{Interp:RN}) and (\ref{Interp:GN}), the desired estimate of the theorem follows.\QED

\begin{remark}\rm

If $J\in L^1_2(\R^n)\cap L^{1+\varepsilon_0}_0(\R^n)$ for some $\varepsilon_0>0$, then by taking $b=0,Q=\infty,q=1$ in (\ref{Est:IntqQ}), we get the pointwise time decay of solution to $\P_tu=\int_{\R^n}J(x-y)(u(y,t)-u(x,t))dx$ to be
\[
\|u(t)\|_{L^\infty}\leq C\lng t\rng^{-n/2}\max\{\|u_0\|_{L^1},\|u_0\|_{L^\infty}\}\quad\forall\,t>0,
\]
for all $u_0\in L^1(\R^n)\cap L^\infty(\R^n)$. 
In \cite{ChasseigneEtal06}, this asymptotic behavior was proved, as the special case $\gamma=2$, assuming 
the Fourier transform of $J$ satisfies 
$\widehat{J}(\xi)=1-A|\xi|^\gamma+o(|\xi|^\gamma)$ as $\xi\to0$ ($0<\gamma\leq2$), and $u_0,\widehat{u}_0\in L^1(\R^n)$.

\end{remark}

\section{Proof of Theorem \ref{Thm:Blow-up}}\label{Sec:Blowup}

To prove the blow-up results for solutions of (\ref{Eqn:main}), 
we begin with the following lemma which is nothing new but a restatement of Theorem \ref{Thm:AsympPsi}.

\begin{lemma}\label{Lem:Blowup1}

Assume $\delta_0\geq2$, $J\geq0$, $J\in L^1_{\delta_0}(\R^n)$, and $|b|\leq\delta_0-2$. Let $R\geq2$. Then there is a constant $d>0$ independent of $R$ such that
\begin{align}\label{Exp:phiR}
\phi_R(x)=\left(1+\frac{|x|^2}{R}\right)^{-b/2},
\end{align}
satisfies
\begin{align}\label{Tmp:LemBU1}
\int_{\R^n}J(x-y)\phi_R(y)dy-\phi_R\geq-\frac{d}{R}\phi_R.
\end{align}

\end{lemma}

We will also need the following lemma.

\begin{lemma}\label{Lem:Blowup2}

Let $f:(0,\infty)\to[0,\infty)$ be a nontrivial differentiable function satisfying
\[
f'(t)\geq-\L f(t)+\mu f(t)^p\quad(t>0),
\]
where $\L\geq0,\mu>0$, $p>1$ are constants. Then $f$ blows up in finite time, i.e.\ $f(t)\to\infty$ as $t\to T_0^-$ for some $T_0<\infty$, provided there exists $t_0\geq0$ such that
\begin{align}\label{BlowupCond}
f(t_0)>\left(\L/\mu\right)^{1/(p-1)}.
\end{align}

\end{lemma}

\noindent\textbf{Proof.} 
We take $t_0$ such that $f(t_0)>0$. 
Let $t\geq t_0$. 
By integration, then $f$ satisfies
\begin{align*}
&f(t)\geq\frac{e^{-\L t}}{\Delta_\L(t)^{1/(p-1)}}\quad(t\geq t_0),\quad\mbox{where}\\
&\Delta_\L(t)
=\begin{cases}
\displaystyle
f(t_0)^{1-p}-(p-1)\mu(t-t_0)&\mbox{if $\L=0$},\\
\vspace{-5pt}\\
\displaystyle
\left[f(t_0)^{1-p}-\frac{\mu}{\L}\right]e^{-(p-1)\L t_0}+\frac{\mu}{\L}e^{-(p-1)\L t}&\mbox{if $\L>0$}.
\end{cases}
\end{align*}

\textbf{Case $\L=0$.} We have
$\Delta_0(t_0)>0$ and $\Delta_0(\infty)=-\infty$. So $f$ always blows up in finite time.\medskip

\textbf{Case $\L>0$.} We have
$\Delta_\L(t_0)>0$ and $\Delta_\L(\infty)=[f(t_0)^{1-p}-(\mu/\L)]e^{-(p-1)\L t_0}$. So $f$ blows up in finite time provided $\Delta_\L(\infty)<0$, i.e.~we can choose $t_0$ such that (\ref{BlowupCond}) holds. \QED

\subsection{Proof of Theorem \ref{Thm:Blow-up} (1)}

Let $\phi_R$ be of the form (\ref{Exp:phiR}) where $b>n$, $R\geq2$ are constants to be specified. Note that $\phi_R\in L^1(\R^n)$ with $\|\phi_R\|_{L^1}=C_{n,b}R^{n/2}$. Multiply (\ref{Eqn:main}) by $\phi_R$, integrate over $\R^n$, and apply (\ref{Tmp:LemBU1}) to get
\begin{align}
\int\phi_R\P_tudx&=\int\phi_R(x)\int\{J(x-y)u(y,t)dy-u(x,t)\}dx+\int\phi_Rau^pdx\nonumber\\
&=\int u(y,t)\int\{J(x-y)\phi_R(x)dx-\phi_R(y)\}dy+\int\phi_Rau^pdx,\nonumber\\
&\geq-\frac{d}{R}\int\phi_Rudx+\int\phi_Rau^pdx\label{BU1:1}.
\end{align}
We have omit $\R^n$ in every integral, for simplicity.\medskip

Let $q=p/(p-1)$. By H\"older's inequality, we have
\begin{align}
&\int\phi_Rudx\leq\left(\int\phi_Rau^pdx\right)^{1/p}\left(\int\phi_Ra^{-q/p}dx\right)^{1/q},\quad\mbox{so}\nonumber\\
&\int\phi_Rau^pdx\geq\left(\int\phi_Rudx\right)^p\left(\int\phi_Ra^{-q/p}dx\right)^{-(p-1)}.\label{BU1:2}
\end{align}
Using (\ref{Hyp:blowup}) and that $p>1$, we estimate
\begin{align}
\int\phi_Ra^{-q/p}dx&\lesssim_a\int\left(1+\frac{|x|^2}{R}\right)^{-b/2}(1+|x|^2)^{-\frac{\S}{2(p-1)}}dx=:I+II+III,\label{BU1:3}
\end{align}
where $I_1,I_2,I_3$, are the respective integrals over $\{|x|<1\}$, $\{1\leq|x|<R^{1/2}\}$, $\{|x|\geq R^{1/2}\}$. Let
\[
f_R(t)=\int\phi_Rudx.
\]
By (\ref{BU1:1})-(\ref{BU1:3}), we obtain
\begin{align}
f_R'(t)&\geq-\frac{d}{R}f_R(t)+C_{a,p}\left(I_1+I_2+I_3\right)^{-(p-1)}f_R(t)^p.\label{Tmp:blowup1}
\end{align}

We estimate $I_1,I_2,I_3$ as follows:
\begin{align}
&I_1\leq C,\label{Tmp:blowup2}\\
&I_2\leq C\int_{1\leq|x|\leq R^{1/2}}|x|^{-\frac{\S}{p-1}}dx
\leq\begin{cases}
\displaystyle
C&\mbox{if $p<1+\frac{\sigma}{n}$},\\
\vspace{-10pt}\\
\displaystyle
CR^{\frac{1}{2}\left(n-\frac{\S}{p-1}\right)}&\mbox{if $p>1+\frac{\sigma}{n}$},\\
\vspace{-10pt}\\
\displaystyle
C\ln R&\mbox{if $p=1+\frac{\sigma}{n}$},
\end{cases}\label{Tmp:blowup3}\\
&I_3\leq R^{n/2}\int_{|z|\geq1}|z|^{-b}(R|z|^2)^{-\frac{\S}{2(p-1)}}dz\leq CR^{\frac{1}{2}\left(n-\frac{\S}{p-1}\right)},\label{Tmp:blowup4}
\end{align}
provided $b>0$ is large enough so that
\begin{align}\label{BU1:4}
b+\frac{\S}{p-1}>n.
\end{align}
We note that all constants $C$'s in the preceding estimates depend only on $n,b,\S,p$ but not on $R$.\medskip

Now we choose and fix $b>n$ such that (\ref{BU1:4}) is true.\medskip

\textbf{Case 1: $1<p<1+\frac{\S}{n}$.} If $\S\leq0$ this case is vacuous, so let us assume $\S>0$. Then we have by (\ref{Tmp:blowup2})--(\ref{Tmp:blowup4}) that $I_1+I_2+I_3\leq C$, hence (\ref{Tmp:blowup1}) becomes
\begin{align*}
f_R'(t)\geq-\frac{d}{R}f_R(t)+C_1f_R(t)^p,
\end{align*}
where $C_1>0$ is independent of $R$. By Lemma \ref{Lem:Blowup2}, $f_R$ blows up in finite time provided
\begin{align}
f_R(t_0)=\int\phi_Ru(t_0)dx>\left(\frac{d}{C_1R}\right)^{1/(p-1)}\quad\mbox{for some $t_0\geq0$}.\label{Tmp:blowup5}
\end{align}
This can be achieved at $t_0=0$ by taking $R$ large enough. In fact, as $R\to\infty$, we have by the monotone convergence theorem that $\int\phi_Ru(0)dx\to\int u_0dx>0$, 
whereas $(d/(C_1R))^{1/(p-1)}\to0$.
\medskip

Now fix $R$ so that (\ref{Tmp:blowup5}) is true at $t_0=0$. By Lemma \ref{Lem:Blowup2}, there is $T_0<\infty$ such that
\[
\int\phi_Ru(t)dx\to\infty\quad\mbox{as $t\to T_0^-$}.
\]
Since $\phi_R\leq1$, we have $\|u(t)\|_{L^1}\to\infty$ as $t\to T_0^-$. On the other hand,
\begin{align*}
\int\phi_Ru(t)dx&\leq\|\phi_R\|_{L^1}\|u(t)\|_{L^\infty}=CR^{n/2}\|u(t)\|_{L^\infty},
\end{align*}
hence $\|u(t)\|_{L^\infty}\to\infty$ as $t\to T_0^-$ as well. Therefore $u$ blows up in finite time in both the $L^1$- and $L^\infty$-norms.\medskip

\textbf{Case 2: $p=1+\frac{\S}{n}$.} In this case $\S>0$. Then we get by (\ref{Tmp:blowup2})--(\ref{Tmp:blowup4}) that $I_1+I_2+I_3\leq C\ln R$, hence (\ref{Tmp:blowup1}) becomes
\begin{align*}
f_R'(t)\geq-\frac{d}{R}f(t)+C_2(\ln R)^{-(p-1)}f_R(t)^p,
\end{align*}
where $C_2>0$ is independent of $R$. The condition for blow-up in finite time is now
\begin{align}\label{Tmp:blowup6}
f_R(t_0)=\int\phi_Ru(t_0)dx>\left(\frac{d}{C_2R}\right)^{1/(p-1)}\ln R.
\end{align}
By the same argument as \textbf{Case 1}, this inequality is true at $t_0=0$ by taking $R$ large, hence the solution $u$ blows up in finite time in both $L^1$- and $L^\infty$-norms.\medskip

\textbf{Case 3: $1+\frac{\S}{n}<p<1+\frac{\S+2}{n}$.} In this case $n>\frac{\S}{p-1}$ and $n<\frac{\S+2}{p-1}$. By (\ref{Tmp:blowup2})--(\ref{Tmp:blowup4}), we have
\[
I_1+I_2+I_3\leq CR^{\frac{1}{2}\left(n-\frac{\S}{p-1}\right)}
\]
and (\ref{Tmp:blowup1}) becomes
\begin{align*}
f_R'(t)\geq-\frac{d}{R}f_R(t)+C_3R^{-\frac{1}{2}\left(n(p-1)-\S\right)}f_R(t)^p,
\end{align*}
where $C_3>0$ is a constant depending only on $n,a,b,\S,p$ but not on $R$. The condition for blow-up in finite time is
\begin{align}\label{Tmp:blowup7}
f_R(t_0)=\int\phi_Ru(t_0)dx
>\left(\frac{d}{C_3}\right)^{1/(p-1)}R^{-\frac{1}{2}\left(\frac{\S+2}{p-1}-n\right)}.
\end{align}
Since $\frac{\S+2}{p-1}-n>0$, we get by the same argument as above that the blow-up condition holds at $t_0=0$ by taking $R$ large. Thus the solution $u$ blows up in both $L^1$- and $L^\infty$-norms.\medskip

Combining the three cases above, we finish the proof of Theorem \ref{Thm:Blow-up} (1).\QED

\subsection{Proof of Theorem \ref{Thm:Blow-up} (2)}

Let $p=1+\frac{\sigma+2}{n}$. We can employ \textbf{Case 3} in the proof of Theorem \ref{Thm:Blow-up} (1) above. By (\ref{Tmp:blowup7}), we get the blow-up in finite time condition as
\begin{align}\label{Tmp:blowupC1}
\int\phi_R(x)u(x,t_0)\,dx>K_0:=\left(\frac{d}{C_3}\right)^{1/(p-1)}\quad\mbox{for some $t_0\geq0$ and $R>0$}.
\end{align}

Assume to get a contradiction that $u\geq0$ is a nontrivial global solution of (\ref{Eqn:main}), which implies in particular that (\ref{Tmp:blowupC1}) is false, i.e.
\begin{align}\label{Tmp:blowupC2}
\int\phi_R(x)u(x,t_0)\,dx\leq K_0\quad\mbox{for all $t_0\geq0$ and $R>0$ large}.
\end{align}
Our aim is to show that $u\equiv0$, which constitutes a contradiction. The contradiction obviously proves the theorem.

\begin{claim}\label{Claim:blowupC}
We have $\displaystyle
\int_0^\infty\int a(x,t)u(x,t)^p\,dxdt\leq K_0<\infty$.

\end{claim}

\noindent\textbf{Proof of Claim.} 
Taking $R\to\infty$ in (\ref{Tmp:blowupC2}) we get 
\[
\int u(x,t_0)\,dx\leq K_0\quad\mbox{for all $t_0\geq0$}.
\]
Integrating (\ref{Eqn:main}) over $\R^n$ and then over time, respectively, we get
\begin{align*}
&\frac{d}{dt}\int u(x,t)\,dx=\int a(x,t)u(x,t)^pdx\\
&\int_0^t\int a(x,\tau)u(x,\tau)^p\,dxd\tau=\int u(x,t)\,dx-\int u_0(x)\,dx\leq K_0,
\end{align*}
this is true for all $t_0\geq0$, so the claim is true.\QED
\medskip

As in \cite{MelianEtal10} (see also  
\cite{MitidieriEtal01}), let $\varphi\in C^\infty_c(\R^n),\psi\in C^\infty_c([0,\infty))$ with $\supp\vp\subset B_1(0)$, $\supp\psi\subset[0,1]$, $0\leq\vp,\psi\leq1$, and 
\[
\vp=1\quad\mbox{on $B_{1/2}(0)$}\quad\mbox{and}\quad\psi=1\quad\mbox{on $[0,1/2]$}.
\]
Fix an arbitrary $t_0\geq0$. Let $R>0$ and  define
\[
\vp_R(x)=\vp\left(\frac{x}{R}\right),\quad\psi_R(t)=\psi\left(\frac{t-t_0}{R^2}\right)\quad(t\geq t_0).
\]
Multiply (\ref{Eqn:main}) by $\vp_R(x)\psi_R(t)$ and integrate over $\R^n\times[t_0,\infty)$ to get
\begin{align}
\int_{t_0}^\infty\int\vp_R\psi_R\P_tu\,dxdt&=\int_{t_0}^\infty\int\vp_R\psi_R(J\ast u-u)\,dxdt+\int_{t_0}^\infty\int\vp_R\psi_Rau^p\,dxdt\nonumber\\
&=\int_{t_0}^\infty\int u\psi_R(J\ast\vp_R-\vp_R)\,dxdt+\int_{t_0}^\infty\int\vp_R\psi_Rau^p\,dxdt,\label{Tmp:blowupC3}
\end{align}
where we have employed the Fubini's theorem.\medskip

Let us consider
\begin{align*}
J\ast\vp_R-\vp_R&=\int_{\R^n}J(y)\left(\vp\left(\frac{x+y}{R}\right)-\vp\left(\frac{x}{R}\right)\right)dy,
\end{align*}
and by Taylor's theorem about $y=0$, we have
\begin{align*}
\vp\left(\frac{x+y}{R}\right)-\vp\left(\frac{x}{R}\right)=\frac{1}{R}y\cdot(\N\vp)\left(\frac{x}{R}\right)+\frac{1}{R^2}\sum_{i,j}y_iy_j(\P_{x_ix_j}\vp)\left(\frac{x}{R}\right)+O(1/R^3).
\end{align*}
The functions $\N\vp,\P_{ij}\vp$ and $O(1/R^3)$ vanish on $|x|<R/2$ and $|x|>R$. The radial symmetry of $J$ implies $\int_{\R^n}J(y)ydy=0$, hence
\begin{align}\label{Tmp:blowupC4}
\left|(J\ast\vp_R-\vp_R)(x)\right|&\leq\frac{C}{R^2}\cdot1_{\{\frac{1}{2}R\leq|x|\leq R\}}.
\end{align}

For the integral on the left of (\ref{Tmp:blowupC3}), we integrate by parts in $t$ to get
\begin{align}\label{Tmp:blowupC5}
\int_{t_0}^\infty\int\vp_R\psi_R\P_tu\,dxdt
&\leq-\int_{t_0}^\infty\int u\vp_R\P_t\psi_R\,dxdt.
\end{align}
Combining (\ref{Tmp:blowupC3})--(\ref{Tmp:blowupC5}) we obtain
\begin{align}
\int_{t_0}^\infty&\int\vp_R\psi_Rau^p\,dxdt\nonumber\\
&\leq\frac{C}{R^2}\int_{t_0}^\infty\int_{\frac{R}{2}\leq|x|\leq R}\psi_Ru\,dxdt-\int_{t_0}^\infty\int u\vp_R\P_t\psi_R\,dxdt\nonumber\\
&\leq\frac{C}{R^2}\left\{\int_{t_0}^{t_0+R^2}\int_{\frac{R}{2}\leq|x|\leq R}u\,dxdt+\int_{t_0+R^2/2}^{t_0+R^2}\int_{|x|\leq R}u\,dxdt\right\},\label{Tmp:blowupC6}
\end{align}
where we have used that $\supp\vp_R\subset B_R$, $\supp\psi_R\subset[t_0,t_0+R^2]$, and $\supp\P_t\psi_R\subset[\frac{1}{2}R^2,R^2]$.\medskip

Using the hypothesis (\ref{Hyp:blowup}), we can estimate
\begin{align*}
\int_{|x|\leq R}a^{-q/p}dx&\leq C\int_{|x|\leq R}(1+|x|^2)^{-\S/(2(p-1))}dx\leq CR^{2/(p-1)}.
\end{align*}
By H\"older's inequality and the preceding inequality, we have
\begin{align}
\int_{t_0}^{t_0+R^2}&\int_{\frac{R}{2}\leq|x|\leq R}u\,dxdt\nonumber\\
&\leq R^{2/q}\left(\int_{t_0}^{t_0+R^2}\left(\int_{\frac{R}{2}\leq|x|\leq R}u\,dx\right)^pdt\right)^{1/p},\quad q:=p/(p-1)\nonumber\\
&\leq R^{2/q}\left(\int_{t_0}^{t_0+R^2}\left(\int_{\frac{R}{2}\leq|x|\leq R}au^p\,dx\right)\left(\int_{\frac{R}{2}\leq|x|\leq R}a^{-q/p}dx\right)^{p/q}dt\right)^{1/p}\nonumber\\
&\leq CR^2\left(\int_{t_0}^{t_0+R^2}\int_{\frac{R}{2}\leq|x|\leq R}au^pdxdt\right)^{1/p}.\label{Tmp:blowupC7}
\end{align}
Similarly, we have
\begin{align}\label{Tmp:blowupC8}
\int_{t_0+R^2/2}^{t_0+R^2}\int_{|x|\leq R}u\,dxdt\leq CR^2\left(\int_{t_0+R^2/2}^{t_0+R^2}\int_{|x|\leq R}au^pdxdt\right)^{1/p}
\end{align}
From (\ref{Tmp:blowupC6})-(\ref{Tmp:blowupC8}) we obtain
\begin{align*}
\int_{t_0}^\infty\int&\vp_R\psi_Rau^p\,dxdt\\
&\leq C\left\{\left(\int_{t_0}^{t_0+R^2}\int_{\frac{R}{2}\leq|x|\leq R}au^p\,dxdt\right)^{1/p}+\left(\int_{t_0+R^2/2}^{t_0+R^2}\int_{|x|\leq R}au^p\,dxdt\right)^{1/p}\right\}
\end{align*}
Taking $R\to\infty$, we get by Claim \ref{Claim:blowupC} that both terms on the right hand side of the preceding estimate converge to zero, hence we conclude that
\begin{align*}
\int_{t_0}^\infty\int au^pdxdt=0.
\end{align*}
It follows that $u\equiv0$ for $t\geq t_0$, but $t_0$ is arbitrary so $u\equiv0$ for all $t$.\QED

\subsection{Proof of Theorem \ref{Thm:Blow-up} (3)}

Let $\phi_R$ be as defined in (\ref{Exp:phiR}). We can follow \textbf{Case 3} in the proof of Theorem \ref{Thm:Blow-up}(1) and get that a positive solution $u$ of (\ref{Eqn:main}) blows up in finite time provided (\ref{Tmp:blowup7}) is true. Since we wish to find how large the initial data $u_0$ is so the solution in this super-critical exponent always blows up, we take $t_0=0$. Now $\frac{\S+2}{p-1}-n<0$, so the right hand side of (\ref{Tmp:blowup7}) is increasing in $R$. To achieve its minimum value, we recall the condition $\eta=R\geq2$ in Theorem \ref{Thm:AsympPsi}. Thus the size of $u_0$ that make the blow-up condition holds at $t_0=0$ is
\[
\int\left(1+\frac{|x|^2}{2}\right)^{-b}u_0(x)dx>m_0:=\left(\frac{d}{C_3}\right)^{1/(p-1)}2^{-\frac{1}{2}(\frac{\S+2}{p-1}-n)}.
\]
Also, recalling that $b>n$ and (\ref{BU1:4}) are requirement for $b$, so we define
\begin{align*}
b_0=\max\left\{n,n-\frac{\S}{p-1}\right\}.
\end{align*}

Thus, for $p>1+\frac{\sigma+2}{n}$, the solution $u$ to (\ref{Eqn:main}) with an initial condition $u_0\geq0$ blows up in finite time provided
\begin{align*}
\int\left(1+\frac{|x|}{2}\right)^{-b}u_0(x)dx>m_0
\end{align*}
for some $b>b_0$.\QED

\section{Proof of Theorem \ref{Thm:Global}}\label{Sec:Global}

We choose constants
\begin{align*}
&b:=\frac{\S}{p-1},\quad\B:=\frac{n-b}{2},\quad\mbox{and}\quad\delta>0\,\,\mbox{small to be specified}.
\end{align*}
Observe that $b\geq0$ and, since $p>1+(\sigma+2)/n$, we have $(p-1)\B>1$. In our investigation of global solutions for (\ref{Eqn:main}), we consider the following Banach spaces
\begin{align}
&\CZ:=C(\R^n)\cap L^1_b(\R^n)\cap L^\infty_b(\R^n),\,\,\mbox{and}\nonumber\\
&\CZ_{gl}(\delta):=\left\{v\in\CZ_{gl}:=C\left([0,\infty);\CZ\right):v\geq0,\|v\|_{\CZ_{gl}}\leq\delta\right\},\,\,\mbox{with the norm}\label{Norm:Global}\\
&\|v\|_{\CZ_{gl}}:=\sup_{t\geq0}\left\{\|v(t)\|_{L^1}+\lng t\rng^{-\frac{b}{2}}\|v(t)\|_{L^1_b}+\lng t\rng^{\frac{n}{2}}\|v(t)\|_{L^\infty}+\lng t\rng^{\B}\|v(t)\|_{L^\infty_b}\right\}.\nonumber
\end{align}
We will often use the fact that
\begin{align}\label{vinZgl}
v\in\CZ_{gl}(\delta)\quad\Rightarrow\quad
\begin{cases}
\displaystyle
\|v(t)\|_{L^1_b}\leq\delta\lng t\rng^{\frac{b}{2}},\quad
\|v(t)\|_{L^1}\leq\delta,\\
\vspace{-10pt}\\
\displaystyle
\|v(t)\|_{L^\infty_b}\leq\delta\lng t\rng^{-\B},\quad
\|v(t)\|_{L^\infty}\leq\delta\lng t\rng^{-\frac{n}{2}}.
\end{cases}
\end{align}
By Theorem \ref{Thm:Interp}, the following estimates are true:
\begin{align}\label{Glo:Interp}
\begin{cases}
Q=q:\quad\|\CG(t)f\|_{L^q_b}\lesssim\lng t\rng^{\frac{b}{2}}\|f\|_{L^q}+\|f\|_{L^q_b},\\
\vspace{-10pt}\\
\displaystyle
\,\,\,\mbox{weight${}=1$}:\quad\|\CG(t)\|_{L^q}\lesssim\|f\|_{L^q},\\
\vspace{-10pt}\\
\displaystyle
q=1,Q=\infty:\quad\|\CG(t)f\|_{L^\infty_b}\lesssim\lng t\rng^{-\B}\|f\|_{L^1}+\lng t\rng^{-\frac{n}{2}}\|f\|_{L^1_b}+e^{-\frac{t}{2}}\|f\|_{L^\infty_b},\\
\vspace{-10pt}\\
\displaystyle
\,\,\,\mbox{weight${}=1$}:\quad\|\CG(t)f\|_{L^\infty}\lesssim\lng t\rng^{-\frac{n}{2}}\|f\|_{L^1}+e^{-\frac{t}{2}}\|f\|_{L^\infty}.
\end{cases}
\end{align}

Define the operator
\[
\CM v(t)=\CG(t)u_0+\int_0^t\CG(t-\tau)\left\{a(\tau)v(\tau)^p\right\}d\tau\quad(t\geq0),
\]
for $v\in\CZ_{gl}(\delta)$. Throughout the following discussion let 
\[
v,w\in\CZ_{gl}(\delta),\quad V=\CM v,\quad W=\CM w.
\]
Clearly $V(t),W(t)\geq0$ for all $t\geq0$.\medskip

\textbf{Step 1.} We show that $
\CM:\CZ_{gl}(\delta)\to\CZ_{gl}(\delta)$, provided $u_0\geq0$ and $\delta>0$ are small enough. 
We use (\ref{vinZgl}) and (\ref{Glo:Interp}) to estimate several norms of $V$.\medskip

\textbf{$L^1_b$-norm of $V(t)$.} Set $(Q,q)=(1,1)$ in (\ref{Glo:Interp}) to get
\begin{align*}
\|V(t)\|_{L^1_b}&\lesssim\lng t\rng^{b/2}\|u_0\|_{L^1}+\|u_0\|_{L^1_b}+\int_0^t\left\{\lng t-\tau\rng^{b/2}\|a(\tau)v(\tau)^p\|_{L^1}+\|a(\tau)v(\tau)^p\|_{L^1_b}\right\}d\tau.
\end{align*}
We estimate
\begin{align}
&\|a(\tau)v(\tau)^p\|_{L^1}
\lesssim\|v(\tau)\|_{L^\infty_b}^{p-1}\|v(\tau)\|_{L^1}\lesssim\delta^p\lng\tau\rng^{-(p-1)\B},\label{Tmp:Glob1}\\
&\|a(\tau)v(\tau)^p\|_{L^1_b}
\lesssim\|v(\tau)\|_{L^\infty_b}^{p-1}\|v(\tau)\|_{L^1_b}\lesssim\delta^p\lng\tau\rng^{-(p-1)\B+\frac{b}{2}},\label{Tmp:Glob2}
\end{align}
where we have used that $|a(x,t)|\leq C\lng x\rng^\sigma$ and (\ref{vinZgl}). 
Thus 
\begin{align}
\lng t\rng^{-b/2}\|V(t)\|_{L^1_b}&\lesssim\|u_0\|_{L^1}+\lng t\rng^{-b/2}\|u_0\|_{L^1_b}+\delta^p\int_0^t\lng\tau\rng^{-(p-1)\B}d\tau\nonumber\\
&\lesssim\|u_0\|_{L^1}+\|u_0\|_{L^1_b}+\delta^p,\label{Tmp:Glob3}
\end{align}
where we have used that $(p-1)\B>1$.\medskip

\textbf{$L^1$-norm of $V(t)$.} Setting $(Q,q)=(1,1)$ and weight${}=1$ in (\ref{Glo:Interp}) and (\ref{Tmp:Glob1}), we have\begin{align}
\|V(t)\|_{L^1}\lesssim\|u_0\|_{L^1}+\delta^p.\label{Tmp:Glob4}
\end{align}

\textbf{$L^\infty_b$-norm of $V(t)$.} Using (\ref{Glo:Interp}) when $(Q,q)=(\infty,1)$, $(Q,q)=(\infty,\infty)$, we get
\begin{align*}
\|&V(t)\|_{L^\infty_b}\\
&\lesssim\lng t\rng^{-\B}\|u_0\|_{L^1}+\lng t\rng^{-n/2}\|u_0\|_{L^1_b}+e^{-t/2}\|u_0\|_{L^\infty_b}\\
&\hspace{0.8cm}+\int_0^{t/2}\left\{\lng t-\tau\rng^{-\B}\|a(\tau)v(\tau)^p\|_{L^1}+\lng t-\tau\rng^{-\frac{n}{2}}\|a(\tau)v(\tau)^p\|_{L^1_b}+e^{-\frac{t-\tau}{2}}\|a(\tau)v(\tau)^p\|_{L^\infty_b}\right\}d\tau\\
&\hspace{2cm}+\int_{t/2}^t\left\{\lng t-\tau\rng^{b/2}\|a(\tau)v(\tau)^p\|_{L^\infty}+\|a(\tau)v(\tau)^p\|_{L^\infty_b}\right\}d\tau.
\end{align*}
Next, we estimate
\begin{align}
&\|a(\tau)v(\tau)^p\|_{L^\infty_b}
\lesssim\|v(\tau)\|_{L^\infty_b}^p\lesssim\delta^p\lng\tau\rng^{-p\B},\label{Tmp:Glob5}\\
&\|a(\tau)v(\tau)^p\|_{L^\infty}
\lesssim\|v(\tau)\|_{L^\infty_b}^{p-1}\|v(\tau)\|_{L^\infty}\lesssim\delta^p\lng\tau\rng^{-(p-1)\B-\frac{n}{2}}.\label{Tmp:Glob6}
\end{align}
Employing (\ref{Tmp:Glob1}), (\ref{Tmp:Glob2}), (\ref{Tmp:Glob5}), and (\ref{Tmp:Glob6}), we get
\begin{align}
\lng t\rng^{\B}\|V(t)\|_{L^\infty_b}&\lesssim\|u_0\|_{L^1}+\lng t\rng^{-b/2}\|u_0\|_{L^1_b}+\|u_0\|_{L^\infty_b}\nonumber\\
&\hspace{1.5cm}+\delta^p\int_0^{t/2}\left\{\lng\tau\rng^{-(p-1)\B}+\lng\tau\rng^{-p\B}\right\}d\tau+\delta^p\int_{t/2}^t\lng\tau\rng^{-p\B}d\tau,\nonumber\\
&\lesssim\|u_0\|_{L^1}+\|u_0\|_{L^1_b}+\|u_0\|_{L^\infty_b}+\delta^p,\label{Tmp:Glob7}
\end{align}
where we used that $\lng t-\tau\rng\sim\lng t\rng$ if $0\leq\tau\leq t/2$ and $\lng t-\tau\rng\lesssim\lng\tau\rng$ if $t/2\leq\tau\leq t$. \medskip

\textbf{$L^\infty$-norm of $V(t)$.} Employing (\ref{Glo:Interp}) with $(Q,q,\mbox{weight})=(\infty,1,1)$, $(Q,q,\mbox{weight})=(\infty,\infty,1)$  together with (\ref{Tmp:Glob1}), (\ref{Tmp:Glob2}), (\ref{Tmp:Glob5}), (\ref{Tmp:Glob6}), we have
\begin{align*}
\|V(t)\|_{L^\infty}&\lesssim\lng t\rng^{-\frac{n}{2}}\|u_0\|_{L^1}+e^{-\frac{t}{2}}\|u_0\|_{L^\infty}\\
&\hspace{0.2cm}+\int_0^{t/2}\left\{\lng t-\tau\rng^{-\frac{n}{2}}\|a(\tau)v(\tau)^p\|_{L^1}+e^{-\frac{t-\tau}{2}}\|a(\tau)v(\tau)^p\|_{L^\infty}\right\}d\tau+\int_{t/2}^t\|a(\tau)v(\tau)^p\|_{L^\infty}d\tau\\
&\lesssim\lng t\rng^{-\frac{n}{2}}(\|u_0\|_{L^1}+\|u_0\|_{L^\infty})+\delta^p\int_0^{t/2}\left\{\lng t-\tau\rng^{-\frac{n}{2}}\lng\tau\rng^{-(p-1)\B}+e^{-\frac{t-\tau}{2}}\lng\tau\rng^{-(p-1)\B-\frac{n}{2}}\right\}d\tau\\
&\hspace{2cm}+\delta^p\int_{t/2}^t\lng\tau\rng^{-(p-1)\B-\frac{n}{2}}d\tau,
\end{align*}
hence we get as above that
\begin{align}
\lng t\rng^{n/2}\|V(t)\|_{L^\infty}\lesssim\|u_0\|_{L^1}+\|u_0\|_{L^\infty}+\delta^p.\label{Tmp:Glob8}
\end{align}

From (\ref{Tmp:Glob3}), (\ref{Tmp:Glob4}), (\ref{Tmp:Glob7}), and (\ref{Tmp:Glob8}), we conclude that
\begin{align*}
\|V\|_{\CZ_{gl}}\leq C_1\left(\|u_0\|_{L^1_b}+\|u_0\|_{L^\infty_b}+\delta^p\right),
\end{align*}
where $C_1$ is a constant independent of $\delta,u_0,V$. Choose $\delta>0$ so that
\begin{align}
\delta\leq\frac{1}{(2C_1)^{1/(p-1)}}.\label{Small1}
\end{align}
Then for any $u_0\in\CZ$ sufficiently small so that
\begin{align}
\|u_0\|_{L^1_b}+\|u_0\|_{L^\infty_b}\leq\frac{\delta}{2C_1}\label{Small2}
\end{align}
we obtain that $\|\CM v\|_{\CZ_{gl}}\leq\delta$. 
Therefore $\CM:\CZ_{gl}(\delta)\to\CZ_{gl}(\delta)$.\medskip

\textbf{Step 2.} We show that $\CM$ is a contraction on $\CZ_{gl}(\delta)$. Let $v,w\in\CZ_{gl}(\delta)$ and $V=\CM v,W=\CM w$. We apply the Theorem \ref{Thm:Interp} and Lemma \ref{Lem:Ineq1} to estimate various norms for
\[
V(t)-W(t)=\int_0^t\CG(t-\tau)\{a(\tau)(v(\tau)^p-w(\tau)^p)\}d\tau.
\]
Observe that $\|v(\tau)\|_{L^\infty_b}\vee\|w(\tau)\|_{L^\infty_b}\leq\delta\lng\tau\rng^{-\beta}$.\medskip

\textbf{$L^1_b$-norm of $V(t)-W(t)$.} By (\ref{Glo:Interp}) with $Q=q=1$, we have
\begin{align}
\lng t\rng^{-b/2}&\|V(t)-W(t)\|_{L^1_b}\leq\lng t\rng^{-b/2}\int_0^t\left\|\CG(t-\tau)\{a(\tau)(v(\tau)^p-w(\tau)^p)\}\right\|_{L^1_b}d\tau\nonumber\\
&\lesssim\int_0^t\left\{\|a(\tau)(v(\tau)^p-w(\tau)^p)\|_{L^1}+\lng t\rng^{-b/2}\|a(\tau)(v(\tau)^p-w(\tau)^p\|_{L^1_b}\right\}d\tau\nonumber\\
&\lesssim\int_0^t(\|v(\tau)\|_{L^\infty_b}\vee\|w(\tau)\|_{L^\infty_b})^{p-1}\left[\|v(\tau)-w(\tau)\|_{L^1}+\lng t\rng^{-b/2}\|v(\tau)-w(\tau)\|_{L^1_b}\right]d\tau\nonumber\\
&\lesssim\delta^{p-1}\|v-w\|_{\CZ_{gl}}\int_0^t\lng\tau\rng^{-(p-1)\B}d\tau\nonumber\\
&\lesssim\delta^{p-1}\|v-w\|_{\CZ_{gl}}\quad(\because\,\,(p-1)\beta>1).\label{Con1}
\end{align}

\textbf{$L^1$-norm of $V(t)-W(t)$.} Using (\ref{Glo:Interp}) with $Q=q=1$ and weight${}=1$, we get
\begin{align}
\|V(t)-W(t)\|_{L^1}&\leq\int_0^t\left\|\CG(t-\tau)\{a(\tau)(v(\tau)^p-w(\tau)^p)\}\right\|_{L^1}d\tau\nonumber\\
&\lesssim\int_0^t\|a(\tau)(v(\tau)^p-w(\tau)^p)\|_{L^1}d\tau\nonumber\\
&\lesssim\int_0^t(\|v(\tau)\|_{L^\infty_b}\vee\|w(\tau)\|_{L^\infty_b})^{p-1}\|v(\tau)-w(\tau)\|_{L^1}d\tau\nonumber\\
&\lesssim\delta^{p-1}\|v-w\|_{\CZ_{gl}}\int_0^t\lng\tau\rng^{-(p-1)\B}d\tau\nonumber\\
&\lesssim\delta^{p-1}\|v-w\|_{\CZ_{gl}}.\label{Con2}
\end{align}

\textbf{$L^\infty_b$-norm of $V(t)-W(t)$.} Applying (\ref{Glo:Interp}) with $(Q,q)=(\infty,1),(\infty,\infty)$, we have
\begin{align}
\lng t\rng^{\B}&\|V(t)-W(t)\|_{L^\infty_b}\lesssim\lng t\rng^{\B}\int_0^t\|\CG(t-\tau)\{a(\tau)(v(\tau)^p-w(\tau)^p)\}\|_{L^\infty_b}d\tau\nonumber\\
&\lesssim\lng t\rng^{\B}\int_0^{t/2}\left\{\lng t-\tau\rng^{-\B}\|a(\tau)(v(\tau)^p-w(\tau)^p)\|_{L^1}+\lng t-\tau\rng^{-\frac{n}{2}}\|a(\tau)(v(\tau)^p-w(\tau)^p)\|_{L^1_b}\right.\nonumber\\
&\hspace{6cm}\left.+e^{-\frac{t-\tau}{2}}\|a(\tau)(v(\tau)^p-w(\tau)^p)\|_{L^\infty_b}\right\}d\tau\nonumber\\
&\hspace{1cm}+\lng t\rng^{\B}\int_{t/2}^t\left\{\lng t-\tau\rng^{\frac{b}{2}}\|a(\tau)(v(\tau)^p-w(\tau)^p)\|_{L^\infty}+\|a(\tau)(v(\tau)^p-w(\tau)^p)\|_{L^\infty_b}\right\}d\tau\nonumber\\
&\lesssim\int_0^{t/2}\Big\{(\|v(\tau)\|_{L^\infty_b}\vee\|w(\tau)\|_{L^\infty_b})^{p-1}\|v(\tau)-w(\tau)\|_{L^1}\nonumber\\
&\hspace{2.5cm}+\lng t\rng^{-\frac{b}{2}}(\|v(\tau)\|_{L^\infty_b}\vee\|w(\tau)\|_{L^\infty_b})^{p-1}\|v(\tau)-w(\tau)\|_{L^1_b}\nonumber\\
&\hspace{4cm}+\lng t\rng^\B e^{-\frac{t-\tau}{2}}(\|v(\tau)\|_{L^\infty_b}\vee\|w(\tau)\|_{L^\infty_b})^{p-1}\|v(\tau)-w(\tau)\|_{L^\infty_b}\Big\}d\tau\nonumber\\
&\hspace{1cm}+\lng t\rng^\B\int_{t/2}^t\Big\{\lng t-\tau\rng^{\frac{b}{2}}(\|v(\tau)\|_{L^\infty_b}\vee\|w(\tau)\|_{L^\infty_b})^{p-1}\|v(\tau)-w(\tau)\|_{L^\infty}\nonumber\\
&\hspace{4cm}+(\|v(\tau)\|_{L^\infty_b}\vee\|w(\tau)\|_{L^\infty_b})^{p-1}\|v(\tau)-w(\tau)\|_{L^\infty_b}\Big\}d\tau\nonumber\\
&\lesssim\delta^{p-1}\|v-w\|_{\CZ_{gl}}.\label{Con3}
\end{align}

\textbf{$L^\infty$-norm of $V(t)-W(t)$.} Finally, applying (\ref{Glo:Interp}) we get
\begin{align}
\lng t\rng^{\frac{n}{2}}&\|V(t)-W(t)\|_{L^\infty}\lesssim\lng t\rng^{\frac{n}{2}}\int_0^t\|\CG(t-\tau)\{a(\tau)(v(\tau)^p-w(\tau)^p)\}\|_{L^\infty}d\tau\nonumber\\
&\lesssim\lng t\rng^{\frac{n}{2}}\int_0^{t/2}\left\{\lng t-\tau\rng^{-\frac{n}{2}}\|a(\tau)(v(\tau)^p-w(\tau)^p)\|_{L^1}+e^{-\frac{t-\tau}{2}}\|a(\tau)(v(\tau)^p-w(\tau)^p)\|_{L^\infty}\right\}d\tau\nonumber\\
&\hspace{2cm}+\lng t\rng^{\frac{n}{2}}\int_{t/2}^t\|a(\tau)(v(\tau)^p-w(\tau)^p)\|_{L^\infty}d\tau\nonumber\\
&\lesssim\lng t\rng^{\frac{n}{2}}\int_0^{t/2}\Big\{\lng t-\tau\rng^{-\frac{n}{2}}(\|v(\tau)\|_{L^\infty_b}\vee\|w(\tau)\|_{L^\infty_b})^{p-1}\|v(\tau)-w(\tau)\|_{L^1}\nonumber\\
&\hspace{2cm}+e^{-\frac{t-\tau}{2}}(\|v(\tau)\|_{L^\infty_b}\vee\|w(\tau)\|_{L^\infty_b})^{p-1}\|v(\tau)-w(\tau)\|_{L^\infty}\Big\}d\tau\nonumber\\
&\hspace{3cm}+\lng t\rng^{\frac{n}{2}}\int_{t/2}^t(\|v(\tau)\|_{L^\infty_b}\vee\|w(\tau)\|_{L^\infty_b})^{p-1}\|v(\tau)-w(\tau)\|_{L^\infty}d\tau\nonumber\\
&\lesssim\delta^{p-1}\|v-w\|_{\CZ_{gl}}.\label{Con4}
\end{align}

Combining the estimates (\ref{Con1})-(\ref{Con4}), we obtain
\begin{align*}
\|\CM v-\CM w\|_{\CZ_{gl}}\leq C_2\delta^{p-1}\|v-w\|_{\CZ_{gl}}
\end{align*}
for some constant $C_2>0$ independent of $u_0,v,w$, and $\delta$. Now we choose $\delta>0$ so that
\begin{align}
\delta\leq\frac{1}{(2C_2)^{1/(p-1)}}\label{Small3}
\end{align}
and (\ref{Small1}) is true. Then $\CM:\CZ_{gl}(\delta)\to\CZ_{gl}(\delta)$ is a contraction.\medskip

\textbf{Step 3.} If $u_0\geq0$ satisfies (\ref{Small2}), then there is a unique fixed point for $\CM$ by the Banach contraction mapping theorem. Therefore (\ref{Eqn:main}) admits a global solution.\QED

\begin{remark}\rm

Observe that if $u_0$ is sufficiently small according to the preceding theorem, the global solution exhibits decay (or extinction):
\[
\|u(t)\|_{L^\infty}\leq\frac{C}{(1+t)^{n/2}},\quad\|\lng\cdot\rng^{\frac{\sigma}{p-1}}u(t)\|_{L^\infty}\leq\frac{C}{(1+t)^{\frac{n}{2}-\frac{\sigma}{2(p-1)}}}.
\]
The latter implies in particular that
\[
|u(x,t)|\leq Ct^{-n/2}\left(\frac{|x|}{\sqrt{t}}\right)^{-\frac{\sigma}{p-1}}\quad\mbox{on $\{|x|\geq\sqrt{t}\}$, for $t\geq t_0>0$},
\]
hence the inhomogeneous coefficient $a(x,t)$ affects the solution so that it decays to zero faster in the region $|x|\geq\sqrt{t}$ at the rate of $\sigma/(p-1)$ in the self-similar variable $|x|/\sqrt{t}$.

\end{remark}

\bigskip

\end{document}